\documentclass[a4paper,12pt]{article}       % onecolumn (second format)

\usepackage[T1]{fontenc}
\usepackage[utf8]{inputenc}

\usepackage[margin=2.5cm]{geometry}
\usepackage{graphicx}
\usepackage[dvipsnames,usenames]{color}
\usepackage{amsmath, amssymb, amsfonts, mathrsfs}

\usepackage{amsthm} % proposition
  \theoremstyle{plain}
  \newtheorem{theorem}{Theorem}
  \newtheorem{proposition}{Proposition}
  
  \theoremstyle{remark}
  \newtheorem{remark}{\texttt{Remark}}

\usepackage[
  justification=centerlast,
  singlelinecheck=true,
  font={small,it},
  labelfont=bf
]{caption}[2004/07/16]
\usepackage{booktabs, graphics, subfig, graphicx}
\usepackage{dashundergaps}
\usepackage{epstopdf}
\usepackage{enumerate}
\usepackage{dblfloatfix}
\usepackage{enumitem}
\usepackage{cite}

\DeclareMathOperator{\derivations}{der}
\DeclareMathOperator{\ran}{ran}
\DeclareMathOperator{\End}{End}

\newcommand{\dfnd}{\overset{\text{def}}{=}}
\newcommand{\lie}[1]{ \{#1\} }
\newcommand{\myDer}[2]{ \frac{ d #1( #2) }{d #2} }%\frac{d #1 \(#2\) }{d #2 } }
\newcommand{\bigLie}[1]{ [#1] }
\newcommand{\myVec}[1]{\boldsymbol{\mathit{#1}}}
\newcommand{\coVec}[1]{\overline{\myVec{#1}}}
\newcommand{\crossLie}[1]{ \big[\!\big[ #1 \big]\!\big] }
\newcommand{\id}{\mathbb{I}}

\begin{document}

\title{
  Hamiltonian Perturbation Theory on a Lie Algebra.
  Application to a non-autonomous Symmetric Top.
}

\author{
  Lorenzo Valvo\thanks{\texttt{valvo@axp.mat.uniroma2.it}}\\
  {\small Dipartimeno di Matematica } \\
  {\small Università degli Studi di Roma ``Tor Vergata''}\\
  {\small Via della Ricerca Scientifica 1} \\
  {\small 00133 Roma, Italy }
  \and Michel Vittot \\
  {\small Centre de Physique Théorique } \\
  {\small Aix-Marseille Université \& CNRS } \\
%  {\small Campus de Luminy, Case 907 } \\
  {\small 163 Avenue de Luminy }\\
  {\small 13288 Marseille Cedex 9, France }
}

\maketitle

\begin{abstract}
We propose a perturbation algorithm 
for Hamiltonian systems on a Lie algebra 
$\mathbb{V}$, so that it can be applied 
to non-canonical Hamiltonian systems. 
Given a Hamiltonian system that preserves a 
subalgebra $\mathbb{B}$ of $\mathbb{V}$, when 
we add a perturbation the subalgebra $\mathbb{B}$ 
will no longer be preserved. We show how to 
transform the perturbed dynamical system to 
preserve $\mathbb{B}$ up to terms quadratic 
in the perturbation. We apply this method to 
study the dynamics of a non-autonomous 
symmetric Rigid Body. In this example our 
algebraic transform plays the role of 
Iterative Lemma in the proof of a 
KAM-like statement.
\end{abstract}

A dynamical system on some set $\mathbb{V}$ 
is a \textit{flow}: a one-parameter 
group of mappings associating to a given element
$F\in\mathbb{V}$ (the initial condition) another element 
$F(t)\in\mathbb{V}$, for any value of the parameter $t$. 

A flow on $\mathbb{V}$ is determined by a linear mapping 
$\mathcal{H}$ from $\mathbb{V}$ to itself. However, the flow 
can be rarely computed explicitly. In perturbation 
theory we aim at computing the flow of \mbox{$\mathcal{H}+\mathcal{V}$}, 
where the flow of $\mathcal{H}$ is known, and $\mathcal{V}$ is
another linear mapping from $\mathbb{V}$ to itself. 

In physics, the set $\mathbb{V}$ is often a Lie algebra: 
for instance in 
classical mechanics \cite{arnold_mathematical_1989},
fluid dynamics and plasma physics \cite{morrison_poisson_1982},
quantum mechanics \cite{sakurai_modern_2017}, kinetic
theory \cite{marsden_hamiltonian_1984}, special and general 
relativity \cite{marsden_covariant_1986}. A dynamical 
system set on a Lie algebra is called
a \textit{Hamiltonian system} after W. R. Hamilton, who
first identified this type of structure in classical 
mechanics.

In classical mechanics, the Lie algebra 
$\mathbb{V}$ is the set of functions over 
a symplectic manifold, with the Lie bracket induced
by the symplectic form \cite{arnold_mathematical_1989}. 
On this type of Lie algebra it is possible to 
introduce particular sets of 
coordinates called 
\textit{canonical coordinates}. But for
many Hamiltonian systems 
(like those that we mentioned above) 
canonical coordinates either 
unavailable \cite{morrison_poisson_1982} 
or undesired \cite{littlejohn_hamiltonian_1982}.

At the same time, canonical coordinates
are needed to perform
perturbation theory in classical mechanics.
A dynamical system is determined
through a function called the Hamiltonian and generally 
denoted by $H$. It is called \textit{integrable} 
if it determines a foliation of the symplectic
manifold into invariant tori. Through perturbation 
theory, one tries to find the tori of the 
Hamiltonian $H+V$ (where V is another 
function, that is the perturbation) usually through a 
series expansion around the tori of $H$. 
An efficient and 
elegant approach to perturbation theory 
was proposed by Kolmogorov
in \cite{kolmogorov_preservation_1954}.
His idea was to conjugate
a perturbed Hamiltonian system to a
new one (named the \textit{Kolmogorov
Normal Form} afterwards) which manifestly
preserves an invariant torus.
Many variants of his theorem have been 
proposed (above all by  
Arnold \cite{arnold_proof_1963} and 
Moser\cite{moser_stable_1973},
whence the name KAM theorem) as well as
generalizations to different settings 
(see for instance \cite{bost_tores_1986}, 
\cite{de_la_llave_kam_2005}, \cite{li_persistence_2002}, 
\cite{alishah_tracing_2012}). In fact, today we speak 
more generically about a KAM theory rather 
than ``the'' KAM theorem. Still, all of 
these approaches require canonical coordinates. 

The Kolmogorov Normal Form
is built in two steps (see for instance
\cite{benettin_proof_1984},
\cite{de_la_llave_tutorial_2001}, 
\cite{broer_kam_2004},
\cite{bergounioux_introduction_2016}).
The first one is often called the 
``Iterative Lemma'': it introduces a map from
the perturbed Hamiltonian $H+V$ to a new one, which
preserves the chosen torus up to an error 
of order $V^2$. 
The second step is to build the Kolmogorov 
Normal Form through a
repeated application of the Iterative Lemma
(hence its name). 

In this work we propose an algebraic approach
to perturbation theory that can be applied 
to any Hamiltonian system, as it requires
only the Lie algebraic structure. 
In this approach, already introduced in
\cite{vittot_perturbation_2004}, the unperturbed 
dynamical system $\mathcal{H}$ is 
required to preserve an invariant
subalgebra $\mathbb{B}$ of the whole 
algebra $\mathbb{V}$. 
We show how to build a first order correction
to a perturbed system $\mathcal{H}+\mathcal{V}$,
so that in the new form it preserves 
$\mathbb{B}$ up to a correction quadratic 
in $\mathcal{V}$. 

Then we 
consider the specific case of a non-autonomous 
symmetric Rigid Body; we call this system
the \textit{Throbbing Top}. It is an example of 
a one and a half degrees of freedom system, with 
canonical coordinates. In the generical
algebraic setting, we were not able to 
provide the equivalents of some key 
elements of KAM theory pertinent
to canonical coordinates, like the 
``homological equation'' and
the ``translation of the actions''. We show 
that our algebraic approach leads naturally
to introduce these elements for the
Throbbing Top.

This paper is organised in four sections. 
In section \ref{sec:algebra}
we recall 
some elements from the theory of 
Lie algebras. In section \ref{sec:kam}
we present the algebraic perturbation 
scheme. In section \ref{sec:tt} 
we study the dynamics of a Throbbing Top, proving a
sort of KAM theorem. Here, the results of the
previous section play the role of the Iterative
Lemma. Finally, in 
section \ref{sec:conclusions} we draw conclusions
and give some hints for further developments.

\section{About Lie Algebras}
\label{sec:algebra}

A \textit{Lie algebra} \cite{arnold_mathematical_1989} is a 
vector space $\mathbb{V}$ over a field $\mathbb{K}$ with
a bilinear operation $\{\,,\}$ (the \textsl{bracket})
which is alternating (here $V,W,Z \in \mathbb{V}$)
\begin{equation}
  \{V,W\}\,=\,-\{W,V\}
\end{equation}
and satisfies the Jacobi identity
\begin{equation}
  \{V,\{W,Z\}\}\,+\,\{W,\{Z,V\}\}\,+\{Z,\{V,W\}\}\,=\,0 
\end{equation}

If on $\mathbb{V}$ we define both a bracket
and an associative product
\begin{equation}
  (A\cdot B)\cdot C\,=\,A\cdot( B\cdot C)
\end{equation}
such that the Leibnitz identity holds,
\begin{equation}
  \{A,\,B\cdot C\}\,=\,\{A,B\}\cdot C\,+\,B\cdot \{A,C\}
\end{equation}
then $\mathbb{V}$ is a \textit{Poisson algebra} \cite{marsden_introduction_2013}.

\bigskip

The space of \textit{derivations} of $\mathbb{V}$ is defined by
\begin{equation}
  \derivations\mathbb{V} \dfnd \big\{\, \mathcal{D} \in \End\mathbb{V}\ \text{s.t.}\ 
  \mathcal{D}\{V,W\}\, =\,
  \{\mathcal{D}V, W\}\, +\, \{V, \mathcal{D}W\}\,, \forall\, V,W \in \mathbb{V} \big\} \label{eq:derivations}
\end{equation}
This space is a Lie algebra on its own with the bracket given by the commutator $[\,,]$,
\begin{equation}
  [\mathcal{F},\mathcal{G}] = \mathcal{F}\mathcal{G} - \mathcal{G}\mathcal{F}
  ,\quad \mathcal{F}, \mathcal{G} \in \End\mathbb{V}
\end{equation}

For any element $F \in \mathbb{V}$, we can consider the
mapping ``bracket with $F$''  
\begin{equation}
  \lie{F}\in \derivations\mathbb{V}, \quad \lie{F}\colon G \mapsto \{F,G\}\,,\,\forall G\in\mathbb{V}
\end{equation}
The image of ``bracket with $F$'' is always a 
derivation; any derivation built in this way is called
an \textit{inner} derivation. A derivation which is not inner is
called \textit{outer}.
Analogously, for any $\mathcal{F}\in\derivations\mathbb{V}$ we can 
consider the mapping ``bracket with $\mathcal{F}$''
\begin{equation}
  [\mathcal{F}]\colon \mathcal{G} \mapsto [\mathcal{F},\mathcal{G}] \,,
  \,\forall \mathcal{G}\in\derivations\mathbb{V}
\end{equation}

\bigskip

Given a derivation $\mathcal{H}$ of $\mathbb{V}$, either inner or outer.
we can define a dynamical system on the algebra,
\begin{eqnarray}
  \label{eq:hamSys}
    \left\{
      \begin{aligned}
        &\dot{F} = \mathcal{H} F \\
        &F(0) = F_0
      \end{aligned}
    \right.
\end{eqnarray}
We call it a \textsl{Hamiltonian system}.  
\textsl{Canonical systems} are a particular example of 
Hamiltonian systems, written in terms of an even-dimensional 
set of coordinates $(q_i,p_i)_{i=1}^n$ such that
\begin{equation}
  \{p_i,p_j\} = 0\,,\quad \{q_i,q_j\} = 0\,,\quad \{p_i,q_j\} = \delta_{ij}
\end{equation}

The formal solution of the system \eqref{eq:hamSys} is 
\mbox{$F(t)=e^{ t\mathcal{H} }F_0$}, and 
$e^{t\mathcal{H}}$ is called the \textit{flow} of $\mathcal{H}$
(see equation \eqref{eq:LieSeries} for the definition of the
exponential operator).

To build a \textsl{non-autonomous dynamical system} we start
from an additive group $\mathbb{G}$ (so that the variable 
$t\in\mathbb{G}$ represents time) and consider the space 
\begin{gather}
  \tilde{\mathbb{V}}\dfnd C^{\infty}(\mathbb{G} \mapsto \mathbb{V} ) \cong \mathbb{V} \otimes C^{\infty}(\mathbb{G}\to\mathbb{R}) \nonumber\\
  \tilde{\mathbb{V}} \owns v(\cdot) \colon t \mapsto v(t) \in \mathbb{V}\,,\, \forall t \in \mathbb{K} \nonumber
\end{gather}
We extend the bracket of $\mathbb{V}$ to $\tilde{\mathbb{V}}$ by the rule,
\begin{equation}
  \forall v, w \in \mathbb{V},\quad [v,w](t) \equiv t \mapsto [v(t),w(t)] 
\end{equation}
and so $\tilde{\mathbb{V}}$ inherits the Lie-algebra structure of $\mathbb{V}$. 
The operator $\partial_t \colon \tilde{\mathbb{V}} \to \tilde{\mathbb{V}}$, 
defined by
\begin{equation}
  \partial_t \colon v(t) \mapsto \myDer{v}{t}
\end{equation}
is a derivation of $\tilde{\mathbb{V}}$: infact, by the linearity of 
$\partial_t$ and the bilinearity of $\{\,,\}$, we have
\begin{equation}
  \partial_t \{v,w\}(t) \,=\, \Big\{\myDer{v}{t},w(t)\Big\} + \Big\{v(t),\myDer{w}{t}\Big\},\,\forall v,w \in \mathbb{V}
\end{equation}
A non-autonomous Hamiltonian system on $\tilde{\mathbb{V}}$ is 
given by
\begin{equation}
  \label{eq:dynamicalSystem}
  \dot{F} = \mathcal{H}F + \partial_t F
\end{equation} 
This choice is made for coherence: if $F$ is time independent, then we 
have the same dynamics of $\mathbb{V}$, while \mbox{$\dot{t} = 1$} as
one would naturally expect.

\section{The Algebraic Perturbation Scheme}
\label{sec:kam}

We start from a Hamiltonian system associated 
to a (not necessarily inner) derivation 
$\mathcal{H}$. In classical mechanics we would require this
system to be integrable; here this notion is replaced by the
existence of a subalgebra 
\mbox{$\mathbb{B}\subset\mathbb{V}$} invariant by $\mathcal{H}$. 
In classical mechanics, $\mathbb{B}$ would be an invariant 
torus.

Then we add a perturbation in the form of an 
inner derivation, 
\mbox{$\mathcal{H}\mapsto\mathcal{H}+\lie{V}$}, 
for some $V\in\mathbb{V}$. Then we show how 
to split the perturbation $V$ into two parts:
one preserving $\mathbb{B}$ and another 
(here denoted by $V_*$) quadratic in $V$. 

To build this splitting, we need a so-called
``pseudo-inverse'' of $\mathcal{H}$.
Being $\mathbb{B}$ invariant by $\mathcal{H}$, it will also
contain its kernel, so that we may hope to be able to 
invert $\mathcal{H}$ on the complementary of $\mathbb{B}$.
Let $\mathcal{R}$ be a projector on $\mathbb{B}$, so
that $\mathcal{N}\dfnd 1-\mathcal{R}$ is a projector
on the complementary
of $\mathbb{B}$. Then we call pseudo-inverse
of $\mathcal{H}$ an operator 
\mbox{$\mathcal{G}\colon\mathbb{V}\to\mathbb{V}$} 
satisfying
\begin{equation}
  \label{eq:easy}
  \mathcal{H}\,\mathcal{G}F\,=\,\mathcal{N}F\,,\quad 
    \forall F\in\mathbb{V}
\end{equation}
In the spirit of classical mechanics, we may 
call $\mathcal{G}F$ a ``generating
function'', because we will use it to transform
away (part of) the perturbation. 
Actually, to perform perturbation theory we only 
need $\lie{\mathcal{G}F}$, which is a derivation.
So we will ask directly for a map 
\mbox{$\Gamma\colon\mathbb{V}\to\derivations\mathbb{V}$}
satisfying 
\begin{equation}
  \label{eq:less_easy}
  [\mathcal{H}] (\Gamma F) = \lie{\mathcal{N}F}
\end{equation}
that is, property \ref{P4} of the following 
Proposition \ref{prop:kam}. 
The simpler case of equation \eqref{eq:easy}
is included in the new one \eqref{eq:less_easy}, 
just by setting \mbox{$\Gamma F\,=\,\lie{\mathcal{G}F}$}.
In fact, recalling the definition \eqref{eq:derivations} of derivation,
\begin{equation}
  \mathcal{H}\lie{\mathcal{G}F} W\,=\,\lie{\mathcal{H}(\mathcal{G}F)}W\,+\,
    \lie{\mathcal{G}F}\mathcal{H}W\,,\quad 
  \forall W\in\mathbb{V}
\end{equation}
we find
\begin{equation}
  \lie{\mathcal{H}\,(\mathcal{G}F)}\,=\,
    \mathcal{H}\lie{\mathcal{G}F} \,-
      \lie{\mathcal{G}F}\mathcal{H}\,=\,
        \bigLie{ \mathcal{H} } (\mathcal{G} F) \,=\,\bigLie{ \mathcal{H} } (\Gamma F )
\end{equation} 
so that 
\begin{equation}
  \lie{NF}\,=\,\bigLie{\mathcal{H}}(\Gamma F)\,
    =\,\lie{\mathcal{H}(\mathcal{G}F)} \iff\, 
      \mathcal{H}\,\mathcal{G}F \,=\,NF
\end{equation}
Now we show how our transformation works.

\begin{proposition}
\label{prop:kam}
  
  Let $\mathbb{B}$ be a Lie subalgebra of $\mathbb{V}$.
  Consider $\mathcal{H} \in \derivations\mathbb{V}$ such that  
  \begin{enumerate}[label=(\roman*)]
    \item \label{P1} $\mathcal{H} \mathbb{B} \subseteq \mathbb{B}$
  \end{enumerate}
  Let $V \in \mathbb{V}$ such that\footnote{
    The case $\lie{V}\mathbb{B}\subseteq \mathbb{B}$ is trivial.
  } $\lie{V}\mathbb{B} \not\subseteq \mathbb{B}$.
  
  Assume to have an operator $\mathcal{R} \colon \mathbb{V} \to \mathbb{B}$ 
  and an operator \mbox{$\Gamma \colon \mathbb{V} \to \derivations\mathbb{V}$} 
  that satisfy the properties
  \begin{enumerate}[label=(\roman*)]
    \setcounter{enumi}{1}
    \item \label{P3} $ \mathcal{R}^2 = \mathcal{R} $ 
    \item \label{P4} $ [\mathcal{H}] (\Gamma F) = \lie{\mathcal{N}F}\,,
	    \,\mathcal{N}\dfnd 1-\mathcal{R}\,,\quad\forall F \in \mathbb{V}$ 
  \end{enumerate} 

  Then 
  \begin{eqnarray}
    e^{ \bigLie{\Gamma V} }( \mathcal{H} + \lie{V} ) = \mathcal{H}_{*} 
    + \lie{V_*} \label{eq:mainStep}\\
    \mathcal{H}_{*} \dfnd \mathcal{H} + \lie{\mathcal{R} V} \label{eq:Hstar}
  \end{eqnarray}
  where $V_{*}$ is a series in $V$ of order quadratic or higher.

\end{proposition}

\begin{remark}
  Hypothesis \ref{P3} states that $\mathcal{R}$ is a projector,
  and is equivalent to ask $\mathcal{N}\mathcal{R}=0$
\end{remark}

\begin{proof}
We start by expanding the l.h.s. of equation~\eqref{eq:mainStep},
\begin{equation}
  e^{ \bigLie{\Gamma V} }( \mathcal{H} + \lie{V} )\, =\,  \mathcal{H}\, +\, \bigLie{\Gamma V} \mathcal{H} 
  \,+\, \sum_{l=2}^\infty \frac{ {\bigLie{\Gamma V}}^l }{l!} \mathcal{H}\, +\, \lie{V}\, +\, ( e^{ \bigLie{\Gamma V} } -1 ) \lie{V} 
\end{equation}
By hypothesis~\eqref{P4}, $\bigLie{\Gamma V} \mathcal{H} = -[\mathcal{H}]\Gamma V = -\lie{\mathcal{N}V} $ so that 
\begin{equation}
  \sum_{l=2}^\infty \frac{ { \bigLie{\Gamma V} }^l }{l!} \mathcal{H}\, =\, 
  - \sum_{l=2}^\infty \frac{ { \bigLie{\Gamma V} }^{l-1} }{l!} \bigLie{ \mathcal{H} } \Gamma V \,
  =\,- \sum_{l=2}^\infty \frac{ { \bigLie{\Gamma V} }^{l-1} }{l!} \lie{\mathcal{N}V} 
  \,=\,- \frac{ e^{ \bigLie{\Gamma V} } - 1 - \bigLie{\Gamma V} }{ \bigLie{\Gamma V} } \lie{\mathcal{N} V} 
\end{equation}
the latter expression being formal. Then 
\begin{eqnarray}
  \nonumber
  \begin{aligned}
  e^{ \bigLie{\Gamma V} }( \mathcal{H} + \lie{V} ) 
    &= \mathcal{H} - \lie{\mathcal{N} V} - \frac{ e^{\bigLie{\Gamma V}} - 1 - \bigLie{\Gamma V} }{ \bigLie{\Gamma V} } \lie{\mathcal{N} V}
     + \lie{V} + ( e^{ \bigLie{\Gamma V} } -1 ) \lie{V} = \\
    &= \mathcal{H} + \lie{\mathcal{R} V}  
     + ( e^{ \bigLie{\Gamma V}} -1 ) \lie{V}
     - \frac{ e^{ \bigLie{\Gamma V} } - 1 - \bigLie{\Gamma V} }{ \bigLie{\Gamma V} } \lie{\mathcal{N} V} 
  \end{aligned}
  \label{eq:non_ne_posso_piu}
\end{eqnarray} 
Now consider the following identity in $\mathbb{V}$,
\begin{equation}
  \bigLie{\Gamma V} \lie{F} = \lie{ (\Gamma V)F } 
\end{equation}
which holds because $\Gamma V$ is a derivation, by definition.
In fact, $\forall l\in\mathbb{N}$,
\begin{equation}
  { \bigLie{\Gamma V} }^l\mathcal{F} = { \bigLie{\Gamma V} }^{l-1} \bigLie{\Gamma V} \mathcal{F} %=\\  
    \,=\, { \bigLie{\Gamma V} }^{l-1} \lie{ (\Gamma V) F }\, =\, 
    { \bigLie{\Gamma V} }^{l-2} \lie{ (\Gamma V)^2 F }\,=\, 
    \ldots\, =\,
    \lie{(\Gamma V)^l F} 
\end{equation}
so that
\begin{equation}
    ( e^{ \bigLie{\Gamma V} } -1 ) \lie{V}\, =\, \sum_{l=1}^\infty \frac{ { \bigLie{\Gamma V} }^l}{l!} \lie{V} =\,
    \sum_{l=1}^\infty \frac{ \lie{ (\Gamma V)^l V } }{l!}\,=\, 
    \lie{ \sum_{l=1}^\infty \frac{(\Gamma V)^l V }{l!} }\, =\, \lie{ (e^{\Gamma V} - 1 ) V } 
\end{equation}
One can proceed analogously to prove that 
\begin{equation}
  \frac{ e^{ \bigLie{\Gamma V} } - 1 - \bigLie{\Gamma V} }{ \bigLie{\Gamma V} } \lie{\mathcal{N} V} = 
    \lie{\frac{ e^{\Gamma V} - 1 - \Gamma V }{ \Gamma V } \mathcal{N} V } 
\end{equation} 
Now, if we inject 
\begin{equation}
  V_{*} = ( e^{\Gamma V} -1 ) V
  - \frac{ e^{\Gamma V} - 1 - \Gamma V }{ \Gamma V } \mathcal{N} V 
  \label{eq:vstar}
\end{equation}
into equation \eqref{eq:non_ne_posso_piu} we recover the thesis \eqref{eq:mainStep}.
\end{proof}

As we discussed in the introduction, in KAM theory 
Proposition \ref{prop:kam} would be called an ``Iterative Lemma'', 
because through its iteration one may build a ``good'' Hamiltonian
$\tilde{H}$ which preserves $\mathbb{B}$ exactly. However, 
to call it an ``Iterative Lemma'' one should also show
that, after a first application of the Lie trasform $e^{\Gamma V}$,
we end up with a system that satisfies 
the original hypothesis of the Lemma again. In our case, 
this means to provide two operators $\mathcal{R}_*$ 
and $\Gamma_*$ that satisfy again hypothesis \ref{P3} 
and \ref{P4} with $\mathcal{H}$ replaced by $\mathcal{H}_*$ and
$V$ replaced by $V_*$. Unfortunately we have not figured out a
general formula to build these operators. However, in the 
next section \ref{sec:tt} we will show how to do it in 
a specific example (see in particular Theorem 
\ref{prop:kam_top2}).

\subsection{Quantitative Estimates}
\label{sec:metric}

For any derivation $\mathcal{A}$ the operator
\begin{equation}
  e^{\mathcal{A}} \,\equiv\,\sum_{n=0}^{\infty}\,\frac{\mathcal{A}^n}{n!}
  \label{eq:LieSeries}
\end{equation}
is called a \textit{Lie series}. Such an expression has only a formal meaning, 
unless we introduce a scale of Banach norms\cite{reed_functional_1981} to 
show that the operator $e^{\mathcal{A}}$ is bounded.

A Banach norm is a function $\|\cdot\|\colon\mathbb{V}\to \mathbb{R}_+$
(where $\mathbb{R}_+$ are the positive real numbers) with properties
\begin{eqnarray}
  &\|\,A\,+\,B\,\|\leq\, \|A\| \,+\, \|B\| \label{prop:norm1} \\
  &\|\lambda A\|\,=\,|\lambda|\|A\| \label{prop:norm2} \\
  &\|A\|\,=\,0\implies A\,=\,0 \label{prop:norm3}
\end{eqnarray}
A scale of Banach norms is a family of norms $\{\|\cdot\|_s\}_{s\in\mathbb{I}}$,
where $s$ is called an index and $\mathbb{I}$ is some set, usually the 
positive integers or the positive reals. For an algebra $\mathbb{V}$ with a scale of
Banach norms indexed by $s\in\mathbb{I}$ we introduce the notation
\begin{equation}
  \mathbb{V}_s = \{ f\in\mathbb{V}\quad \text{s.t.}\quad \|f\|_s \leq \infty \}
\end{equation}
and we assume that 
\begin{equation}
  W \in \mathbb{V}_{s_1} \implies W \in \mathbb{V}_{s_2}, \forall s_2 < s_1
\end{equation}

We say that a derivation $\mathcal{D}$ is \textit{bounded with loss} if 
\begin{equation}
  \|\mathcal{D}A\|_{s-\delta} \leq \alpha(\delta) \|A\|_s
\end{equation}
for any $A\in\mathbb{V}_s$, $s,\delta\in\mathbb{I}$. A paradigmatical example,
regularly used in KAM theory \cite{giorgilli_quantitative_1995}, is the following: 
on the complex plane $\mathbb{C}$ we define the sets
\begin{equation}
  \mathbb{B}_r(0)\dfnd\{ z\in\mathbb{C}\,\text{s.t.}\, |z|< r\}
\end{equation}
Then, on the space $\mathcal{C}^{\omega}(\mathbb{C})$ 
we consider the scale of norms (indexed by $r\in\mathbb{R}_+$) 
\begin{equation}
  | f(z) |_r\dfnd \sup_{z\in\mathbb{B}_r} |f(z)|
\end{equation}
The Cauchy inequality states that
\begin{equation}
  \label{eq:cauchy}
  |\partial_z f(0)|\, \leq\, \frac{1}{r}\, |f(z)|_r 
\end{equation}
from which we get the upper bound
\begin{equation}
  |\partial_z f|_{r-\delta}\,\leq \,\frac{1}{\delta}\,|f|_r
\end{equation}
So by loosing a ``layer'' of width $\delta$ of the original
domain, corresponding to a shift $r\mapsto r-\delta$ of the index
of the norm, it was possible to bound from above the derivation
operator on $\mathbb{C}$.

In next Proposition we make the formal manipulations
of Proposition \ref{prop:kam} quantitative by assuming
that the Lie algebra $\mathbb{V}$ is endowed with a 
scale of Banach norms.

\begin{proposition}
\label{prop:kam2}

  Let the Lie algebras $\mathbb{V}$ and $\mathbb{B}$,
  the function $V\in\mathbb{V}$ and the operators 
  $\mathcal{H}$, $\mathcal{R}$, $\mathcal{N}$
  and $\Gamma$ be as in Proposition \ref{prop:kam}.

  Assume that on $\mathbb{V}$ there exists a
  scale of Banach norms $\{\|\quad\|_r\}_{r\in\mathbb{I}}$,
  such that \mbox{$\|V\|_s\, < \infty$} for some $s\in\mathbb{I}$.

  Assume also that $\forall s, d,\delta\in\mathbb{I}, d<\delta, d+\delta< s$
  there exist two functions $\Lambda(d,\delta)$ and $\Xi(\delta)$
  such that
  \begin{enumerate}[label=(\roman*)]
    \item \label{P5} $\| (\Gamma V) F \|_{s-\delta-d}\, \leq\, \Lambda(d,\delta) \, 
	    \|V\|_s\, \| F\|_{s-\delta}\,,\quad F\in\mathbb{V}_{s-\delta}$
    \item \label{P6} $\|\mathcal{N}V\|_{s-\delta}\, <\, \Xi(\delta) \|V\|_s$
    \item \label{eq:epsilon_equation}
            $\epsilon_\delta = \frac{1}{2}\,\sup_{n\in\mathbb{N}} \bigg( \frac{1}{n!}
               \prod_{j=1}^n \Lambda\Big(\,\frac{\delta}{n}\,,\, \frac{(j-1)\delta}{n}\, \Big) 
                         \bigg)^{-1/n}\,\in\,(0,\infty)$
  \end{enumerate}
  Then the operator $e^{ \bigLie{\Gamma V} }$ is well defined, and for any 
  $\mathbb{I} \owns \mu< s/3$ we have
  \begin{equation}
    \label{eq:epsilonStar}
    \|V\|_s \leq \epsilon_\mu\,\implies\,
      \|V_*\|_{s-3\mu} \, \leq\,\kappa\, {\epsilon_\mu}^2
  \end{equation}
  for some real positive constant $\kappa$.

\end{proposition}

\begin{proof}
  We will show that $e^{\Gamma V}$ is bounded with loss
  from $\mathbb{V}_s$ to $\mathbb{V}_{s-\mu}$ (it's easier to 
  study convergence on an algebra rather than on the space of 
  its derivations). Then $e^{\bigLie{\Gamma V}}$ 
  can be computed by the relation,
  \begin{equation}
    \label{eq:relation}
    e^{\bigLie{\mathcal{A}}}\, \mathcal{B} =\, e^{\mathcal{A}}\, \mathcal{B}\,
     \,e^{-\mathcal{A}}, \quad \forall \mathcal{A},\mathcal{B} \in \derivations\mathbb{V}
  \end{equation}
  which is readily proven by using a series 
  expansion on both sides. Indeed we can use the relation 
  \begin{equation}
    \bigLie{\mathcal{A}}^N \mathcal{B} = \sum_{k=0}^N \binom{N}{k} \,\mathcal{A}^k\, \mathcal{B} \,(-\mathcal{A})^{N-k}
  \end{equation}
  to rewrite the l.h.s. of \eqref{eq:relation} as
  \begin{equation}
    \label{eq:relation_1} 
    \sum_{N\geq 0}\, \sum_{k=0}\, \frac{{\mathcal{A}}^k\, \mathcal{B}\, (-\mathcal{A})^{N-k} }{k!\,(N-k)!} 
  \end{equation}
  The r.h.s. of equation \eqref{eq:relation} is
  \begin{equation}
    \sum_{n\geq 0,m\geq 0} \frac{ \mathcal{A}^n}{n!}\, \mathcal{B}\, \frac{(-\mathcal{A})^m}{m!}
  \end{equation}
  and by a change of variable $m\mapsto N-n$ becomes
  \begin{equation}
   \sum_{N\geq 0 }\,\sum_{n=0}^N\, \frac{\mathcal{A}^n\, \mathcal{B}\, (-\mathcal{A})^{N-n} }{n! (N-n)!}
  \end{equation}
  By renaming an index, the above is equal to 
  expression \eqref{eq:relation_1}.
  
  Now consider the expression $(\Gamma V)^n F$, as $n$ varies. For $n=1$ we can 
  apply hypothesis~\ref{P5} with \mbox{$\delta=0$} and \mbox{$d=\mu$} to get 
  \begin{equation}
    \|(\Gamma V) F \|_{s-\mu}\,\leq\,\Lambda(\mu,0)\, \|V\|_s\, \|F\|_s
  \end{equation}
  Now let \mbox{$n\geq 1$} and for any \mbox{$1\leq j \leq n$}, consider the operator
  \begin{equation}
    (\Gamma V)^j\,\colon\,\mathbb{V}_{s - (j-1)\mu/n}\,\to\,\mathbb{V}_{s - j\mu/n}
  \end{equation}
  By applying hypothesis~\ref{P5} with \mbox{$d=\mu/n$} and \mbox{$\delta=(j-1)\mu/n$}
  we get
  \begin{equation}
    \| (\Gamma V)^j F \|_{s-\mu}\,
    \leq\,\Lambda( \tfrac{\mu}{n}, \tfrac{(j-1)\mu}{n} ) \,\|V\|_s\, \|(\Gamma V)^{j-1} F\|_{s-(j-1)\mu/n} \\[1ex]
  \end{equation}
  and, iterating the above $n$ times,
  \begin{equation}
    \big\|(\Gamma V)^n F\big\|_{s-\mu}\,\leq\,
      \prod_{j=1}^n \Lambda(\tfrac{\mu}{n}, \tfrac{(j-1)\mu}{n} ) {\|V\|_s}^n  \| F\|_{s}
  \end{equation}
  We can finally bound $e^{\Gamma V}$ with loss, 
  \begin{multline}
    \nonumber
    \big\|e^{\Gamma V}F\big\|_{s-\mu} 
    \,\leq\sum_{n=0}^\infty\,\frac{1}{n!}\,\big\|(\Gamma V)^n F\big\|_{s-\mu}
    \,\leq\\\leq\,
    \sum_{n=0}^\infty \frac{1}{n!}
      \prod_{j=1}^n \Lambda(\tfrac{\mu}{n}, \tfrac{(j-1)\mu}{n} ) {\|V\|_s}^n  \| F\|_s
    \leq\,\sum_{n=0}^\infty \bigg( \frac{\|V\|_s}{ 2 \epsilon_\mu} \bigg)^n \| F\|_{s}\, \leq \,2\|F\|_s
  \end{multline} 
  where we also used equation \eqref{eq:epsilon_equation} and the 
  hypothesis that $\|V\|_s< \epsilon_\mu$.
  
  To bound the norm of $V_*$ we use a similar technique,
  \begin{eqnarray}
    \begin{aligned}
      \nonumber
      \| V_* \|_{s-3\mu}\, &=\, \Big\| \sum_{l\geq 1 } \frac{(\Gamma V)^l}{l!} V - \sum_{l\geq 2} \frac{(\Gamma V)^{l-1} }{l!} \mathcal{N}V \Big\|_{s-3\mu}\, \leq \\[1ex]
      \nonumber
        &\leq\, \sum_{l \geq 1} \frac{1}{l!} \Big\| (\Gamma V)^l V \,-\,\frac{(\Gamma V)^{l-1} }{l+1} \mathcal{N}V \Big\|_{s-3\mu} \,\leq \\[1ex]
      \nonumber
        &\leq\, \sum_{l\geq 0 } \frac{1}{l!} \Big\| (\Gamma V)^l \, \frac{1}{l+1} \Big( (\Gamma V)V \,+\,\frac{1}{l+1} (\Gamma V) \mathcal{N}V \Big) \Big\|_{s-3\mu} \,\leq \\[1ex] 
      \nonumber
        &\leq\, \,\sum_{n=0}^\infty \bigg( \frac{\|V\|_s}{ 2 \epsilon_{\mu} } \bigg)^n \Big( \big\|(\Gamma V)V \big\|_{s-2\mu}\, +\, \| (\Gamma V)\mathcal{N}V \|_{s-2\mu} \Big)\, \leq \\[1ex]
        &\leq\, 2\,\big( \Lambda( 2\mu,0 ) \|V\|_s^2\, +\, \Lambda( \mu,\mu ) \|V\|_s 
          \underbrace{ \|\mathcal{N}V\|_{s-\mu} }_{\leq\,\Xi(\mu)\|V\|_s } \big)\,\leq\,\kappa\, \epsilon_{\mu/3}^2
    \end{aligned}
  \end{eqnarray}
  where we used hypothesis \ref{P5}, \ref{P6} and also that, for any positive integer $l$, $1/(l+1)<1$,
  $1/(l+2)<1$. So we proved equation~\eqref{eq:epsilonStar}.
\end{proof}

The two hypothesis that are usually assumed in KAM
theory, besides the analyticity of the involved functions, are 
the Diophantine condition for the frequency on the torus and 
the non-degeneracy of the Hamiltonian. All of these assumptions
are ``hidden'' into hypothesis \ref{P5} on the existence and 
boundedness of $\Gamma$. Indeed, in section \ref{sec:norms}, we 
use both of them to prove the boundedness of the 
operator $\Gamma$ specific to the Throbbing Top.

\section{The dynamics of a symmetric and periodic Throbbing Top}
\label{sec:tt}

In what follows, we will use the names ``Rigid Body'' or ``Top'' 
as synonyms. However, we prefer the term Top: as we are
considering a non-autonomous system, it is unlikely to be 
``rigid''.

\subsection{Basic facts on the (static) Top}

The space $\mathbb{R}^3$ is a Lie algebra with the bracket $\crossLie{\quad}$
(the vector product). It is also a metric space; we denote by an overbar Euclidean
transposition. As a consequence of the Lie-Poisson theorem (see for
instance \cite{marsden_introduction_2013}), the set
\begin{equation}
  \label{eq:lieAlgebraR3}
  \mathbb{V}_\text{Top}\,\dfnd\,\mathcal{C}^{\infty}( \overline{\mathbb{R}^3} \to \mathbb{R} )
\end{equation}
is a Poisson algebra with bracket
\begin{equation}
  \label{eq:bracket}
  \{ F , G \}\big(\coVec{M}\big) = \coVec{M} \crossLie{ \partial_{\coVec{M}} F} \partial_{\coVec{M}} G
  ,\quad\forall F,G\in\mathbb{V}_\text{Top}
\end{equation}
The operator $\partial_{\coVec{M}}$ on $\mathbb{V}$ is defined by
\begin{equation}
  \coVec{N} \partial_{\coVec{M}} f = \lim_{\eta \to 0} \frac{ f( \coVec{M} + \eta \coVec{N} ) - f( \coVec{M} ) }{ \eta }
\end{equation}
and it takes elements of $\mathbb{V}_\text{Top}$ into elements of $\mathbb{R}^3$. This is evident from
the definition: when we act on $\partial_{\coVec{M}} f$ with 
an element \mbox{$\coVec{N}\in{\mathbb{R}^3}^*$}, we get a scalar.

If we consider as Hamiltonian the function
\begin{equation}
  \label{eq:hamiltonian_rb}
  E\, =\, \frac{1}{2}\, \coVec{M} \mathsf{L} \myVec{M},\quad
    \mathsf{L} = \begin{pmatrix}
                I_1^{-1} & 0 & 0 \\
                0 & I_2^{-1} & 0 \\
                0 & 0 & I_3^{-1}
              \end{pmatrix}
\end{equation}
then by \mbox{$\dot{M} = \lie{E} M$} we recover the Euler-Poinsot equation
for the Rigid Body.

The matrix $\mathsf{L}$ is called the ``tensor of inertia'', and it 
encodes the properties of the Top (its shape, mass distribution \dots). 
This matrix is symmetric, 
so it has three real eigenvalues $\{1/I_i\}_{i=1}^3$, 
that are the inverse of the ``moments of inertia''. And 
these eigenvalues are always positive. In general an 
ordering like \mbox{$I_1 > I_2 > I_3$} or the opposite 
\mbox{$I_1 < I_2 < I_3$} is assumed. The special cases 
\mbox{$I_1=I_2=I_3$} and \mbox{$I_1=I_2$} (or \mbox{$I_2=I_3$}) 
are respectively known as the \textsl{spherical Top}, and 
as the \textsl{symmetric Top}. 

The function
\begin{equation}
  \rho^2\dfnd M_1^2\,+\,M_2^2\,+\,M_3^2 \label{eq:Casimir} 
\end{equation}
represents the (square) modulus of $\myVec{M}$ ans has the property
\mbox{$\lie{\rho}F=0$}, for any \mbox{$F\in\mathbb{V}_\text{Top}$}; we call
it a \textit{Casimir element} \cite{marsden_introduction_2013}. 
A Casimir element is constant under the flow determined
by any Hamiltonian; in fact, it is a property of the algebra, not of the 
flow. As a consequence, the dynamics of a Top 
takes place in a two-dimensional space: a sphere of radius $\rho$. 
It is possible to show that, given a Hamiltonian system 
on a Poisson algebra, after quotienting away the Casimir 
elements, we get a canonical system. And a two dimensional, 
autonomous canonical system is integrable\footnote{
  In the context of sympletic mechanics, a dynamical system of dimension
  $2n$ is called integrable if it has $n$ quantities in involution
  (i.e. having zero bracket) among themselves and with the Hamiltonian. 
  As an obvious consequence, a canonical Hamiltonian system
  is always integrable for $n=1$, which is the case of the Top.
}, and so is the case for the Top. But a non-autonomous system is
no longer integrable, even in two dimensions.

\begin{figure}
    \centering
    \includegraphics[height=0.3\textheight]{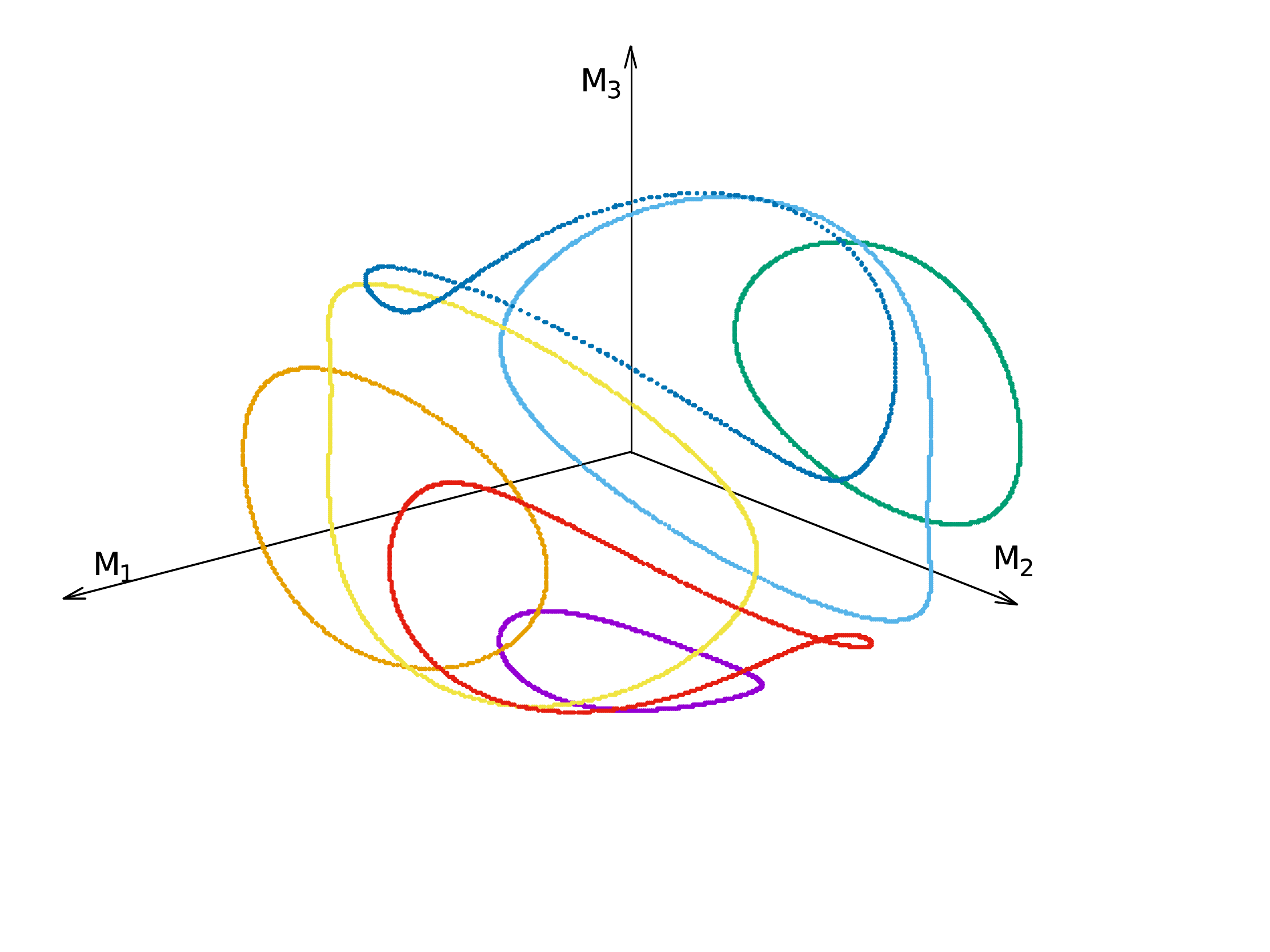}
    \caption{
      A few trajectories of a static Top 
      with moments of inertia \mbox{$I_1=1$}, \mbox{$I_2=2$}, 
      \mbox{$I_3=3$}. 
      These trajectories were generated by 
      a code emplying a Runge-Kutta 4th order integration
      scheme and step $h=0.001$. The initial data
      were randomly generated with the unique constraint 
      of having all the same value of $\rho=2$. Conservation
      of $\rho$ and of the energy
      was achieved up to numerical precision.
    }
    \label{fig:ex_rigid_body}
\end{figure}

In the static case, the energy $E$ and the Casimir 
$\rho$ determine two surfaces in $\mathbb{R}^3$, 
a sphere and an ellipsoid, so that the intersections of the two objects give
trajectories of the Top. As a consequence, there exists a set of 
accessible values for the energy: given $\rho$ and the moments of inertia, the 
system will have a solution only for 
\begin{equation}
  \rho^2/(2I_3)\, \leq\, E\, \leq\, \rho^2/(2I_1)
  \label{eq:enBound}
\end{equation}
(if $I_1 > I_3$). In figure~\ref{fig:ex_rigid_body} we plot a few trajetories
for a Rigid Body with moments \mbox{$I_i = i$}.

\subsection{The Throbbing Top}

The mathematical description of a non-autonomous periodic Top, according
to section \ref{sec:algebra}, is set on the algebra
\begin{equation}
  \label{eq:algebraTT}
  \mathbb{V}_{TT}\,\dfnd\, \mathcal{C}^\infty (\,
     \mathbb{T}\, \to\,\mathbb{V}_\text{Top}\,)\, \owns\, f = f\big(\coVec{M},t\big)
\end{equation}
again with the bracket \eqref{eq:bracket}. As the time variable doesn't 
enter in the bracket, $\rho$ is still a Casimir. This means that the 
energy, even if it is fluctuating, has to respect the bound
\eqref{eq:enBound}. The phase space, that in the static case was the
sphere $\mathbb{S}_2$, becomes \mbox{$\mathbb{S}_2\times \mathbb{T}$}.

We will assume that the unperturbed Hamiltonian is still given 
by \eqref{eq:hamiltonian_rb}. We are interested in perturbations of type
\begin{equation}
  \label{eq:top_perturbation}
  V = \frac{1}{2}\, \coVec{M} \mathsf{A}(t) \myVec{M}
\end{equation}
where $\mathsf{A}(t)$ is a $3\times 3$ diagonal matrix with time dependent coefficients.
Physically, this will represent a Top for which the moments of inertia are 
changing in time.  

\begin{figure}
    \centering
    \includegraphics[height=0.3\textheight]{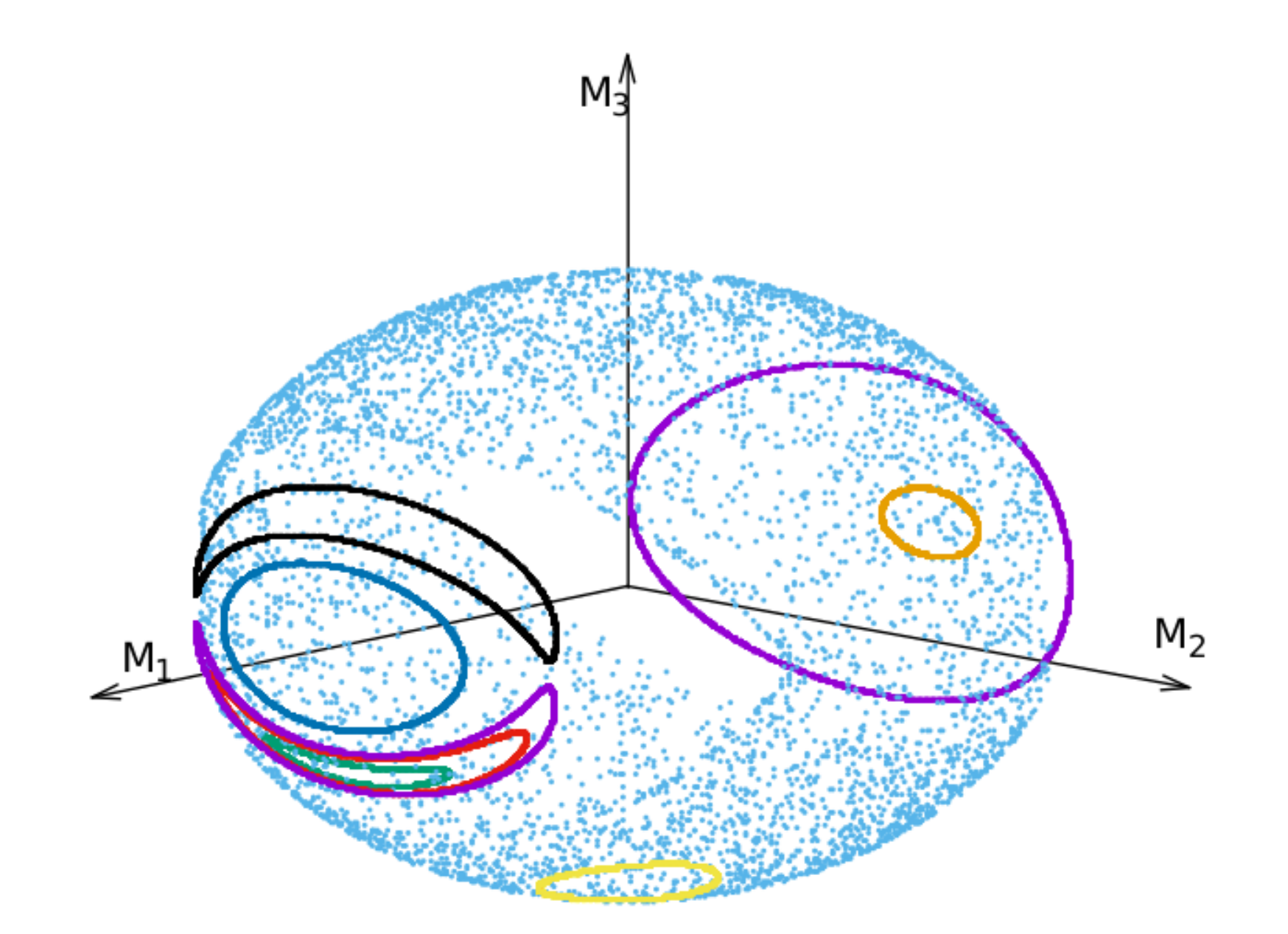}
    \caption{
      A few trajectories of a Throbbing Top 
      with $I_1=1,\,I_2=2/( 1 + 0.2 \epsilon\cos(t)),\,I_3=3 $.    
      The initial data were generated as in \ref{fig:ex_rigid_body},
      and the numerical method was the same as well. We checked 
      also the conservation of $\rho$ along these trajectories.
    }
    \label{fig:ex_t-d-top}
\end{figure}

The new dynamical system is 
\begin{equation}
  \label{eq:throbbingTop}
    \dot{F}\,=\,\mathcal{H}F\,+\,\lie{V}F\,,\quad \mathcal{H}\,=\,\lie{E}\,+\,\partial_t
\end{equation}
For instance, for $F=M_i,\, i=1,2,3$ and $V$ given by 
\eqref{eq:top_perturbation} with
\begin{equation}
  \mathsf{A}(t) = \begin{pmatrix}
            0 & 0 & 0 \\
            0 & \epsilon \cos(\nu t) & 0 \\
            0 & 0 & 0 
         \end{pmatrix} 
  \label{eq:pert1}
\end{equation}
we are describing a Top with
\mbox{$I_2 \,= \, I_2^\emptyset/( 1 + I_2^\emptyset \epsilon\cos(\nu t))$}
being $I_2^\emptyset$ the static value of $I_2$. 
In figure~\ref{fig:ex_t-d-top} 
we plot some trajectories of this dynamical system.
We observe the typical features of dynamical systems with 
cohexistence of order and chaos. The separatrices (the lines joining
the hyperbolic equilibria $M_1=0,M_3=0$) disappear, and 
are replaced by orbits spanning a two-dimensional area. Around 
the elliptic equilibrium points (of coordinates respectively
$M_1=0,M_2=0$ and $M_2=0,M_3=0$), some of the original trajectories
are only deformed, some others are lost and replaced by a set 
of new equilibrium points; some of the new equilibrium points are 
elliptic, and new closed orbits appear around them. 

\subsection{The symmetric case}
\label{sec:symm}

By definition a Top is symmetric if two moments of inertia are
equal; here we fix $I_1=I_2\equiv I_\perp$. 
In this case the solutions of motion % \eqref{eq:EulerPoinsot}
are uniform rotations around the third axis
($M_3$ is constant in time). 

For a static symmetric top it is useful \cite{deprit_free_1967}, 
\cite{gurfil_serretandoyer_2007} to introduce the coordinates 
$(\rho,X,\theta)$, $X\in (-1,1)$ and $\theta\in[0,2\pi)$, 
defined by
\begin{equation}
  \label{eq:coordinates}
  \{M_1, M_2, M_3\} \mapsto \{\rho,X,\theta\}\,\colon%\\[1ex]
  \left\{
  \begin{aligned}
  M_1 &= \rho \sqrt{1-X^2} \cos( \theta ) \\[0.8ex]
  M_2 &= \rho \sqrt{1-X^2} \sin( \theta ) \\[0.8ex]
  M_3 &= \rho X
  \end{aligned}
  \right.
\end{equation}
The bracket \eqref{eq:bracket} restricted to $\mathbb{V}_\text{symm}$ 
becomes\footnote{
  by abuse of notation, we use the same symbol $\lie{\quad}$ as before
}
\begin{equation}
  \label{eq:bracket2}
  \lie{F}G = \frac{1}{\rho} \big( \partial_X\! F\, \partial_\theta G\, -\, \partial_\theta F\, \partial_X G \big)
\end{equation}
The new bracket contains no derivatives in $\rho$, consistent with the 
definition of a Casimir\footnote{And, from this moment on, we won't write $\rho$
anymore among the coordinates.}. The Hamiltonian~\eqref{eq:hamiltonian_rb} becomes
\begin{equation}
  \label{eq:symmHam}
  E_\text{symm}\,=\,\frac{\rho^2}{2}\bigg( \frac{1-X^2}{ I_\perp} + \frac{ X^2}{I_3}\bigg)\, \equiv\,
    \frac{\rho^2}{2}\Delta\,X^2\,+\,\frac{\rho^2}{2  I_\perp}
\end{equation}
where we have set $\Delta=\frac{1}{I_3}-\frac{1}{I_\perp}$.

So we see that $X$ and $\theta$ behave like action-angle coordinates. 

The coordinates $X$ and $\theta$ don't cover
the whole sphere, as the north and south poles are excluded. However, in 
the stationary case the poles are elliptic equilibria, so they are not 
very interesting for the dynamics. In the non-autonomous case the energy 
is still subject to the bound~\eqref{eq:enBound}. If the bound
is strengthened to strict inequalities then the dynamics will
never reach the poles.

So, we restrict the algebra $\mathbb{V}_{TT}$ (defined in \eqref{eq:algebraTT}) 
to the subalgebra of functions $f(X,\theta,t)$ analytic in $(X,\theta,t)$
and which respect the bounds \eqref{eq:enBound} with strict inequality.
The restriction to analytic functions is needed to introduce a scale of 
Banach norms, as will be discussed in subsection \ref{sec:norms}. 

We start by making a further change of coordinates (sometimes called
localization of $X$),
\begin{equation} 
  X\,=\,x_0 + x\,\implies\,\partial_X \mapsto \partial_x
\end{equation}
where $x_0\in(-1,1)$ is fixed and $x$ is sufficiently small so that
$x_0+x\in(-1,1)$. This change of variables is simply a translation
and it doesn't affect the algebraic and metric properties that
we introduced up to now. Functions in $\mathbb{V}_{\text{symm}}$ can be
equivalently written as $f(x,\theta,t)$. Let us also define
\begin{equation}
  \label{eq:Q}
  Q\dfnd \partial^2_{xx}H|_{x=0}
\end{equation} 
so that, for instance,
 $\lie{E_\text{symm}} \,=\,\rho x_0 \Delta \partial_\theta + \lie{\tfrac{1}{2}Qx^2}$

Now we will show that all the hypothesis of Lemma \ref{prop:kam}
are satisfied for the symmetric Throbbing Top, that is, by system
\eqref{eq:throbbingTop} on the algebra $\mathbb{V}_\text{symm}$ 
with $E=E_\text{symm}$.

\bigskip

\textbf{1.} First we look for a subalgebra $\mathbb{B}$ of $\mathbb{V}_\text{symm}$,
invariant by $\mathcal{H}$. Led again by analogy with classical mechanics, we choose
\begin{equation}
  \mathbb{B} = \big\{ F(\rho,x,\theta,t) \in \mathbb{V}\quad  \text{s.t.} \quad F(\rho, 0, \theta, t) = 0, \, \partial_x F(\rho,0,\theta,t) = 0 \big\}
  \label{eq:B}
\end{equation}
By definition we have $F \in \mathbb{B} \iff \mathcal{P}_{\geq 2} F_2 = F_2$
(see Table \ref{tab:op} for the definition of $\mathcal{P}_{\geq 2}$), but
neither $\lie{E_{\text{symm}}}$ nor $\partial_t$ can decrease the degree 
in $x$ of a polynomial. So $\mathcal{H}\mathbb{B}\subseteq\mathbb{B}$.

\begin{table}[t]
  \[
  \begin{array}{c}
  \toprule
  \oint \colon f(x,\theta,t) \mapsto \int_0^{2\pi} d\theta \int_0^{2\pi} dt f(x,\theta,t) \equiv f_{0,0}(x)  
  \\[2ex]
  \chi \dfnd 1 - \oint \colon f \mapsto \sum_{l,m\in\mathbb{Z}_0} f_{l,m}(x) e^{im\theta + ilt},%\\ &
    \mathbb{Z}_0 \dfnd \mathbb{Z}\backslash\{0\} 
  \\[2ex]
  \mathcal{P}_{k} \colon \sum_{n\geq 0 } a_n x^n \mapsto a_k x^k \,\, k \in \mathbb{N},\,a_n \in \mathbb{R}\ \forall n\in \mathbb{N} 
  \\[2ex]
  \mathcal{P}_{\geq k} \colon \sum_{n\geq 0 } a_n x^n \mapsto \sum_{n \geq k} a_n x^n
  \\[2ex]  
  \mathcal{R}_s \dfnd \oint \mathcal{P}_0 + \mathcal{P}_{\geq 2} 
  \\[2ex]
  \mathcal{N}_s \dfnd \chi \mathcal{P}_0 + \mathcal{P}_{1} \equiv 1 - \mathcal{R}_s 
  \\[2ex] 
  \mathcal{G}_s \colon f \mapsto \sum_{l,m \in \mathbb{Z}_0 } \frac{ -i \mathcal{P}_{\leq 1} f_{l,m}(x) }{ \rho x_0 \Delta m + l } e^{ i m \theta + i l t }  
  \\[2ex] 
  \mathcal{A} \dfnd \Big( \oint Q\Big )^{-1} \oint \mathcal{P}_0 ( \partial_x - Q \partial_\theta \mathcal{G}_s) 
  \\[2ex]
    \mathcal{K}\, \dfnd\, \rho^2 x_0 \Delta \mathcal{A}\,+\, \lie{ \tfrac{1}{2}Qx^2 } \, x \mathcal{G}_s\, \big( \mathcal{P}_0 \partial_x  
     -\, Q \mathcal{A}  - Q \partial_\theta \mathcal{G}_s \mathcal{P}_0 \big) \\ 
  \bottomrule
  \end{array}
  \]
  \caption{
    Here we group some of the 
    operators defined on $\mathbb{V}$ and needed for the KAM algorithm.
    We are using 
    the Fourier representation~\eqref{eq:repr1}, the $1$s are to be intended 
    as identity operators, $H$ is the Hamiltonian, 
    and $Q$ has been defined in equation~\eqref{eq:Q}. 
    \label{tab:op}
  }
\end{table}

\bigskip
\textbf{2.} As a second step we build the projector $\mathcal{R}$
(and thus $\mathcal{N} = \id -\mathcal{R}$).
We choose 
\begin{equation}
  \label{eq:RN}
  \mathcal{R} \dfnd \mathcal{R}_s - \mathcal{K},
  \quad \mathcal{N} \dfnd \mathcal{N}_s + \mathcal{K}
\end{equation}
where $\mathcal{R}_s, \mathcal{N}_s$ and $\mathcal{K}$ are defined in
table~\ref{tab:op}. It's evident that $\mathcal{R}_s$ takes values in 
$\mathbb{B}$, and then 
$\mathcal{R} = \mathcal{R}_s - \mathcal{K} \equiv \mathcal{R}_s( 1-\mathcal{K} )$.
In point (\textbf{4}) we show that $\mathcal{N}\mathcal{R}=0$ so 
we can conclude that they are both projectors.

\bigskip
\textbf{3.} The third step is to build the operator $\Gamma$.
As we discussed at the beginning of section \ref{sec:kam},
it would be simpler to compute
\mbox{$\mathcal{G}\colon \mathbb{V}\to\mathbb{V}$} and then
$\Gamma F = \lie{\mathcal{G}F}$.
The operator $\mathcal{G}_s$ from table \ref{tab:op} satisfies
\[
  \big(\rho x_0 \Delta \partial_\theta\,+\,\partial_t)\, \mathcal{G}_s
    F\,=\,\mathcal{N}_s F\,,\quad \forall F\in\mathbb{V}
\]
This equation is called ``homological
equation'' in classical mechanics. Unfortunately,
the term in parenthesis above doesn't correspond
to our $\mathcal{H}$, which needs a more complicated
pseudo-inverse. 

\bigskip
\textbf{4.}
Still making reference 
to table \ref{tab:op}, consider the following operator:
\begin{equation}
  \mathcal{G}\, =\, \mathcal{G}_s\, +\, \rho\, \theta \mathcal{A}\, -\, 
    x\mathcal{G}_s Q ( \mathcal{A} + \partial_\theta \mathcal{G}_s \mathcal{P}_0 )
  \label{eq:G}
\end{equation}
It acts on elements of $\mathbb{V}$, but it doesn't take values in $\mathbb{V}$,
because the function $\theta$ doesn't
belong\footnote{
  This is commonly seen in KAM theory: given a phase space with
  action-angle coordinates $(\myVec{\varphi},\myVec{A})$ the translation of the action
  of a quantity $\myVec{\xi}$ is generated by a function of type 
  $\chi = \coVec{\varphi}\myVec{\xi}$, which doesn't belong to the algebra
  of functions $f(\myVec{A},\myVec{\varphi})$ over the phase space. Indeed, the latter
  functions are periodic in $\myVec{\varphi}$, while this is not the case for $\chi$.

} to $\mathbb{V}$. But we can formally compute 
\begin{equation}
  \Gamma f = \lie{\mathcal{G} f} = \lie{\mathcal{G}_s f} - \rho^{-1} \mathcal{A}f \partial_x -
    \lie{ x\mathcal{G}_s Q ( \mathcal{A} + \partial_\theta \mathcal{G}_s \mathcal{P}_0 ) f }
  \label{eq:Gamma}
\end{equation} 
so $\Gamma$ goes from $\mathbb{V}$ to $\derivations\mathbb{V}$.
So we proceed to check equation \eqref{eq:easy} that in this context reads
\begin{equation}
  ( \rho x_0 \Delta \partial_\theta + \partial_t + \lie{\tfrac{1}{2}Qx^2} )\,
    ( \mathcal{G}_s + \rho \theta \mathcal{A} - x \mathcal{G}_s Q ( \mathcal{A} + \partial_\theta \mathcal{G}_s \mathcal{P}_0 ) ) f = \mathcal{N}f \\
\end{equation}
We use $\lie{ \tfrac{1}{2}Q x^2 } = xQ\partial_\theta - \tfrac{1}{2}x^2(\partial_\theta Q)\partial_x$
and we get
\begin{multline}
  \nonumber
  \chi\mathcal{P}_{\leq 1} f + \lie{ \tfrac{1}{2}Qx^2 } \mathcal{G}_s f 
    -x \chi \mathcal{P}_{\leq 1} Q ( \mathcal{A} + \partial_\theta \mathcal{G}_s \mathcal{P}_0 ) f +
  \underset{\bullet}{\underline{ x Q \mathcal{A} f }} 
  \\ \nonumber  
    +\underset{\sim}{\underline{ \rho x_0 \Delta \mathcal{A} f }}
      + \underset{\square}{\underline{ \lie{ \tfrac{1}{2} Q x^2 } x \mathcal{G}_s Q ( \mathcal{A} + \partial_\theta \mathcal{G}_s \mathcal{P}_0 ) f}} =
  %\\ \nonumber 
  (\chi \mathcal{P}_0 + \mathcal{P}_1 ) f + \underset{\sim}{\underline{ \rho x_0 \Delta \mathcal{A} f}} +
  \\ \nonumber 
    - \underset{\square}{\underline{ \lie{ \tfrac{1}{2} Q x^2 } x \mathcal{G}_s Q \mathcal{A} f }} + 
  \underset{\bullet}{\underline{ x Q \mathcal{A} f}} + \lie{ \tfrac{1}{2}Qx^2 } x \mathcal{G}_s \mathcal{P}_0 \partial_x f -
    \underset{\square}{\underline{ \lie{ \tfrac{1}{2}Qx^2 } x \mathcal{G}_s \mathcal{P}_0 Q \partial_\theta \mathcal{G}_s f }}
\end{multline}
All the terms underlined in the same way cancel among themselves, and we are left with
\begin{equation}
  \lie{ \tfrac{1}{2}Qx^2 } \mathcal{G}_s f + xQ\mathcal{A}f-x(\chi Q) \mathcal{A} f-x\chi Q \partial_{\theta} \mathcal{G}_s \mathcal{P}_0 f\,=
    \,\oint\mathcal{P}_1 f + \lie{ \tfrac{1}{2}Qx^2 } x \mathcal{G}_s \mathcal{P}_0 \partial_x f 
\end{equation}
Now we observe that \mbox{ $x\mathcal{P}_0 \partial_x f = \mathcal{P}_1 f$ } so that 
there is a partial cancellation among the first and the latter term in 
the above equation, and there remains
\begin{equation}
  \lie{ \tfrac{1}{2}Qx^2 } \mathcal{G}_s \mathcal{P}_0 f\, + \,x \big(\oint Q \big) \mathcal{A} f\, -\,x\chi Q \partial_{\theta} \mathcal{G}_s \mathcal{P}_0 f
    =\,\oint\mathcal{P}_1 f 
\end{equation}
Then we insert the explicit expressions of \mbox{$\lie{\tfrac{1}{2}Qx^2}$} and 
that of $\mathcal{A}$ as it can be found in table~\ref{tab:op},
\begin{equation}
  \underset{\diamond}{ \underline{ xQ\partial_\theta \mathcal{G}_s \mathcal{P}_0 f } }\, 
  +\, \underset{\triangle}{\underline{ x \oint \mathcal{P}_0 \partial_x f }}\, 
  -\, \underset{\diamond}{\underline{ x\oint Q \partial_{\theta} \mathcal{G}_s \mathcal{P}_0 f}} \,
  -\, \underset{\diamond}{\underline{ x\chi Q \partial_{\theta} \mathcal{G}_s \mathcal{P}_0 f }}
  \,= \,\underset{\triangle}{\underline{\oint\mathcal{P}_1 f }}
  %\nonumber
\end{equation}
Again we underlined in the same way all the terms that cancels out.
We conclude that equation~\eqref{eq:easy} is satisfied.

\bigskip

\textbf{5.} Here we show that $\mathcal{G}\mathcal{R} = 0$, so that 
$\mathcal{H}\mathcal{G}\mathcal{R} = \mathcal{NR} = 0$. We start by
writing explicitly
\begin{equation}
  \nonumber
  \mathcal{G}\mathcal{R} = 
    \Big(\mathcal{G}_s +\theta \mathcal{A}- x\mathcal{G}_s Q 
    ( \mathcal{A}  + \partial_\theta \mathcal{G}_s \mathcal{P}_0 ) \Big)
      %\times \\ \times
        \Big(\mathcal{R}_s - \rho \Delta x_0 \mathcal{A} -
        \lie{\tfrac{1}{2}Qx^2} x \mathcal{G}_s (\mathcal{P}_0\partial_x -
        Q \mathcal{A} - Q \partial_\theta \mathcal{G}_s \mathcal{P}_0 ) \Big)
\end{equation}
Next we observe that, by applying the following equalities
\begin{eqnarray}
    \nonumber
    \mathcal{G}_s \mathcal{R}_s \propto \chi \mathcal{P}_{\leq 1} (\oint \mathcal{P}_0 + \mathcal{P}_{\geq 2} ) ) = 0 \\
    \nonumber
    \mathcal{A R}_s f = \Big( \oint Q \Big)^{-1} \Big( \oint \mathcal{P}_0 \partial_x \mathcal{P}_{\geq 2} - \oint Q \partial_\theta \mathcal{G}_s \mathcal{R}_s \Big) = 0\\
    \nonumber
    \mathcal{A A} \propto \mathcal{A} \oint \mathcal{P}_0 = 0\\
    \nonumber
    \mathcal{G}_s \mathcal{A} \propto \mathcal{G}_s \oint \mathcal{P}_0 = 0
\end{eqnarray}
many terms cancel, and we are left with
\begin{equation}
  \mathcal{G}\mathcal{ R} = 
    \Big[ - \mathcal{G}_s - \theta \mathcal{A} + x\mathcal{G}_s Q 
    ( \mathcal{A} + \partial_\theta \mathcal{G}_s \mathcal{P}_0 ) \Big]\,
     \underbrace{ \lie{\tfrac{1}{2}Qx^2}x\mathcal{G}_s(\mathcal{P_0}\partial_x - Q \mathcal{A} - Q \partial_\theta \mathcal{G}_s \mathcal{P}_0 ) }_{\in\ran(\mathcal{R}_s) } = 0
\end{equation}
We can conclude that Proposition \ref{prop:kam} can be
applied to the symmetric periodic Throbbing Top.

\subsection{A scale of Banach norms for $\mathbb{V}_{\text{symm}}$}
\label{sec:norms}

Functions in $\mathbb{V}_\text{symm}$ are analytic and thus admit the 
Fourier representation
\begin{equation}
  \label{eq:repr1}
  F(x,\theta,t)\,=\,\sum_{l,m\in\mathbb{Z}} F_{l,m}(x)\,e^{ilt+im\theta} 
\end{equation}
Analyticity allows to build a complex extension $(-1,1)$, the domain
of $X$ (and of $x$). Let $\mathbb{B}_r(X)\subseteq\mathbb{C}$ be a ball
in the complex plane, of radius\footnote{
  we denote by $\mathbb{R}_+$ the set of positive reals.
} $r\in\mathbb{R}_+$ centered at $X$. The radius $r$ has to be sufficiently 
small so that $|X\pm r|<1$.  Then we define the set 
\begin{equation}
  \mathbb{A}_r \dfnd \bigcup_{X\in(-1,1)} \mathbb{B}_r(X)
\end{equation}
The algebra $\mathbb{V}_\text{symm}$ is a subalgebra of 
\begin{equation}
  \mathbb{V}_r\,\dfnd\, \mathcal{C}^\infty( \mathbb{A}_r \otimes \mathbb{T}^2 )
\end{equation}
for any $r$. Moreover, we restrict to the subset of analytic functions, so that
each space $\mathbb{V}_r$ is endowed with the Banach norm
\begin{equation}
  \label{eq:normSymm}
  \big\| f \big\|_{r} \dfnd \sum_{l,m \in \mathbb{Z}} \big| f_{l,m} \big|_r\, e^{r(|l| + |m|)}\,,
  \quad \big| f_{l,m} \big|_r \dfnd \sup_{X \in \mathbb{A}_r} |f(X)|
\end{equation}
So we have a scale of Banach norms $\{\big\|\quad\big\|_r \}_{r\in\mathbb{R}_+}$ and 
a scale of Banach spaces $\{\mathbb{V}_r\}_{r\in\mathbb{R}_+}$. Some properties of 
these norms are collected in the following Proposition (the proof can be easily
reconstructed by adapting the proof of Lemma 1 of \cite{giorgilli_quantitative_1995}).

\begin{proposition}
\label{prop:normProp}

  Consider the Lie algebra $\mathbb{V}_\text{symm}$ 
  with the scale of Banach norms \eqref{eq:normSymm}.
  Let $r,\delta,d \in\mathbb{R}_+$ with $d+\delta < r$. 
  Let also $W \in \mathbb{V}_{r},\,Z\in\mathbb{V}_{r - \delta}$ .
  Then
  \begin{eqnarray}
    \| \partial_X W \|_{r - d}\,\leq\,\frac{1}{d} \| W \|_{r} \label{eq:formula1} \\
    \| \partial_\theta W \|_{r - d}\,\leq\,\frac{1}{ed} \| W \|_{r} \label{eq:formula2} \\
    \| \lie{W} Z \|_{r - d - \delta}\,\leq\,\frac{2}{\rho e d (d + \delta)} \| W \|_{r} \| Z \|_{r-\delta} \label{eq:formula3} 
  \end{eqnarray}
\end{proposition}

Instead in the appendix \ref{app:proof2} we prove the following 

\begin{proposition}
\label{prop:bounds_rb}

  Consider the Lie algebra $\mathbb{V}_\text{symm}$ 
  with the scale of Banach norms \eqref{eq:normSymm}.
  Let \mbox{$V,Q \in \mathbb{V}_r$} for some \mbox{$r\in\mathbb{R}_+$}.
  Define two operators $\mathcal{R},\mathcal{N}$ by~\eqref{eq:RN}
  and an operator $\Gamma$ by \eqref{eq:Gamma}.
  Assume there exist real numbers {$\gamma,\rho,\Delta >0$},
  \mbox{$\tau > 1$},
  \mbox{$0 < q < 1$} and
  \mbox{$-1 < x_0 < 1$}  
  such that:
  \begin{enumerate}[label=(\arabic*)]
    \item \label{1} $x_0,\rho,\Delta,\gamma$ and $\tau$ satisfy
    %\begin{equation}
    %  \label{eq:diophantine}
     $\, \big|\,\rho \Delta x_0 m\, +\, l\,\big|\, \geq\, \gamma\,\big(|l| + |m| \big)^{-\tau},\quad\forall l,m\in\mathbb{Z}_0 $;
    %\end{equation}
   
    \item \label{2} $|Q_{00}|\, \geq\, q$; 
    \item \label{3} $\|Q\|_r\, \leq\, q^{-1}$;
  \end{enumerate}
  
  Then \mbox{$\forall d,\delta\in\mathbb{R}_+, d+\delta< r$}
  and \mbox{$\,\forall W\in\mathbb{V}_r\,,\, \forall Z\in\mathbb{V}_{r-\delta}$},
  the following inequalities hold
  \begin{eqnarray}
    \label{eq:bound}
    \|(\Gamma W)Z \|_{r -\delta-d}\,\leq\,\frac{C \, \|W\|_{r}\,\|Z\|_{r-\delta}  }{q^3 d(d+\delta)^{2\tau+2} } \\
    \label{eq:boundN}
    \| \mathcal{N} W \|_{r-\delta} \leq \frac{ \tilde{C} \, \|W\|_{r} }{ q^3\, \delta^{2\tau+3}} \\
    \label{eq:boundR}
    \| \mathcal{R} W \|_{r-\delta} \leq \frac{ \tilde{C} \,\|W\|_{r} }{q^3\,\delta^{2\tau+3}}
  \end{eqnarray}
  where $C$ and $\tilde{C}$ are constants 
  depending on $\tau,\gamma,e,q,\rho,\delta,d$. 

\end{proposition} 

The inequalities \eqref{eq:bound} and \eqref{eq:boundN},
are respectively of type \ref{P5} and \ref{P6} of
Proposition \ref{prop:kam2}, with
\begin{equation}
  \Lambda(d,\delta)\,=\,\frac{ C % (\tau,\gamma,r,\delta,d,q) 
    }{ q^3\,d\,(d+\delta)^{2\tau+2} }\,, \qquad
  \Xi(\delta)\,=\,\frac{ \tilde{C} }{ q^3\,\delta^{2\tau+3} }
\end{equation} 
We see that some extra hypothesis on $Q$ and on the product $\rho\Delta x_0$
are required. In particular condition~\ref{1} of Proposition \ref{prop:bounds_rb} 
is usually called the ``Diophantine condition''. Instead hypothesis \ref{2} 
goes generally under the name of ``non-degeneracy condition''. Finally,
condition \ref{eq:epsilon_equation} of Proposition \ref{prop:kam2} defines
the parameter $\epsilon$, that in this case equals 
\begin{equation}
  \begin{aligned}
    \nonumber
    \epsilon_\mu\,&=\,\frac{1}{2}\sup_{n\in\mathbb{N}}\bigg( \frac{1}{n!}\,\prod_{j=1}^n \Lambda\Big( \tfrac{\mu}{n},\tfrac{(j-1)\mu}{n} \Big) \bigg)^{-1/n}\, =\,
    \frac{1}{2}\sup_{n\in\mathbb{N}}\bigg( \frac{1}{n!}\,\prod_{j=1}^n \frac{C}{q^3 \,\tfrac{\mu}{n}\,\big(\tfrac{j\mu}{n}\big)^{2\tau+2} }  \bigg)^{-1/n}\,=\\[1.2ex]
    &=\,\frac{ q^3 \mu^{2\tau+3} }{ 2\, C }\,\sup_{ n\in\mathbb{N} }\Big( \frac{n^n}{n!} \Big)^{2\tau+3}\, 
    \label{eq:epsilonSymmetricTop}
    =\,\frac{ q^3 \mu^{2\tau+3} }{ 2\, C }\,\sup_{ n\in\mathbb{N} }\Big( e^{1-1/n} \Big)^{2\tau+3} 
    =\,\frac{ q^3 \mu^{2\tau+3} }{ 2\, C }
  \end{aligned}
\end{equation}
So, if $q,\mu,\tau$ and $C$ are chosen so that $0<\epsilon_\mu<\infty$,
we can conclude that also Proposition \ref{prop:kam2} applies to the
symmetric and periodic Throbbing Top.

\subsection{A KAM theorem for the Symmetric Throbbing Top}
%\label{sec:kam_sptt}

Now we prove a KAM theorem for the symmetric Throbbing Top by 
iteratively applying Proposition \ref{prop:kam}.

\begin{theorem}
  \label{prop:kam_top2}
  Consider the dynamical system \eqref{eq:throbbingTop}
  on the algebra $\mathbb{V}_\text{symm}$ and with 
  $E=E_\text{symm}$ (Throbbing Top). Define $Q$ as in
  equation \eqref{eq:Q} and $\rho,\Delta,\gamma,\tau,q\in\mathbb{R}_+$ 
  as in Proposition~\ref{prop:bounds_rb}.
  Then there exist $\epsilon_0,r\in\mathbb{R}_+$ such that
  if $\|V\|_r \leq \epsilon$, the for a large class of 
  initial data, the trajectories of the Throbbing Top 
  can be mapped to trajectories of a 
  static Top $\mathcal{H}_\infty$. 
\end{theorem}

\begin{proof}

We have shown in the previous section that Proposition \ref{prop:kam} 
can be applied to the symmetric and periodic Throbbing Top, by 
choosing $\mathbb{B}$ as in \eqref{eq:B}, $\mathcal{R}$ as in 
\eqref{eq:RN} and $\Gamma$ as in \eqref{eq:Gamma}. Thus by formula 
\eqref{eq:mainStep} we can map \mbox{$\mathcal{H}+\lie{V}$} into 
\mbox{$\mathcal{H}_*+\lie{V_*}$}. This is possible, in particular, 
if $X(t=0)=x_0$ satisfy the Diophantine condition \ref{1} of 
Proposition \ref{prop:bounds_rb}.. We have to choose a loss 
$\mu_0\in\mathbb{R}_+$ such that $\mu_0<r/3$, so that
$\|V_*\|_{r-3\mu}<\kappa\,\epsilon_0^2$. 

Now we want to show that Proposition \ref{prop:kam} can
be applied to $\mathcal{H}_*+\lie{V_*}\equiv\,\mathcal{H}+\lie{RV}+\lie{V_*}$, 
so with the same values of $\rho,\Delta$ and $x_0$ but   
with the replacements
\begin{equation}
  Q\,\to\,Q^* \dfnd Q + \partial^2_{xx}(\mathcal{R}V)
  \,,\quad\quad r \to r_* \equiv r -\mu\,,\quad\quad
  V\to V_*
\end{equation}
We need to verify that there exist a new 
constant $q_*$ for which hypothesis \ref{2} and \ref{3} of 
Proposition~\ref{prop:bounds_rb} are again satisfied. By using 
inequality~\eqref{eq:boundR} of Proposition~\ref{prop:bounds_rb},
\begin{equation}
  \|Q^* - Q\|_{r-\mu}\, =\, \|\partial^2_{xx}\mathcal{R}V \|_{r-\mu}\, 
    \leq\, \frac{2}{r^2} \| \mathcal{R}V\|_r \,\leq\, \frac{2\,\epsilon\, \tilde{C}}{r^2 q^3 \mu^{2\tau+3} }
\end{equation}
where we used equations \eqref{eq:formula1}, \eqref{eq:boundR}, and 
the Cauchy inequality \eqref{eq:cauchy}. In the same way
\begin{equation}
  |Q^*_{0,0} - Q_{0,0}|\, =\, \Big|\oint \partial^2_{xx}\mathcal{R}V\Big|\, \leq\, \frac{2}{r^2} \Big|\oint \mathcal{R}V \Big|_r\, \leq\, \frac{2\,\epsilon\,\tilde{C}}{ r^2 q^3 \mu^{2\tau+3} }
\end{equation}
So we have
\begin{equation}
  \|Q^*\|_{r-\mu} \,\leq\, \|Q\|_{r_0-\mu}\, +\, \|Q^* - Q\|_{r_0-\mu} \,%\leq \\
    \leq\, \frac{1}{q}\, +\, \frac{\epsilon \tilde{C}}{r^2 q^3 \mu^{2\tau+3} }\, \leq\, \frac{1}{ q - \frac{\epsilon\tilde{C}}{r^2 q^3 \mu^{2\tau+3} } }\, \equiv\, \frac{1}{q_*}
\end{equation}
where the last inequality holds as long as 
\begin{equation}
  0\leq\frac{\epsilon\tilde{C}}{ r^2 q^3 \mu^{2\tau+} }\leq q \leq 1
\end{equation} 
This condition is to be confronted with formula \eqref{eq:epsilonSymmetricTop}; 
they are compatible if
\begin{equation}
  \label{eq:convergenceCondition}
  q\, \,\geq \, \tilde{C}/(C\,r^2)
\end{equation}
At the same time,
\mbox{$|Q_{0,0}|\,\leq\, |Q_{0,0}-Q^*_{0,0}|\, +\, |Q^*_{0,0}|$}
so 
\begin{equation}
  |Q^*_{0,0}|\, \geq\, |Q_{0,0}|-|Q^*_{0,0} - Q_{0,0}| \,\geq\, q\, -\, \frac{\epsilon\tilde{C}}{ q^3 \mu^{ 2\tau + 5 } }\, \equiv\, q_* 
\end{equation}
So Proposition~\ref{prop:kam} can be applied to $\mathcal{H}_* +\lie{V_*}$. 
We may build a sequence of dynamical systems by
\begin{eqnarray}
    e^{\bigLie{\Gamma_i V_i}}(\mathcal{H}_i + \lie{V_i}) = \mathcal{H}_{i+1} + \lie{V_{i+1}} \\[1ex]
    V_0 \equiv V\,,\quad Q^0 \equiv Q\,,\quad \mathcal{H}_i = x_0 \rho \Delta \partial_\theta + \partial_t + \lie{ \tfrac{1}{2} Q^i x^2} \\[1ex]
    \Gamma_i f = \lie{\mathcal{G} f} = \lie{\mathcal{G}_s f} - \rho^{-1} \mathcal{A}_i f \partial_x -
      \lie{ x\mathcal{G}_s Q^i ( \mathcal{A}_i + \partial_\theta \mathcal{G}_s \mathcal{P}_0 ) f } \\[1ex]
    \mathcal{A}_i = \Big(\oint Q^i\Big)^{-1} \oint \mathcal{P}_0( \partial_x - Q^i \partial_\theta \mathcal{G}_s \mathcal{P}_0 ).
\end{eqnarray}
The sequence converges to the static Top
\begin{equation}
  \mathcal{H}_\infty\,=\,\mathcal{H}\,+\,\sum_{i=0}^\infty \mathcal{R}_i V_i
\end{equation}
To show that the sequence exists, we need
three sequences $\{\epsilon_i,\mu_i,q_i\}_{i\in\mathbb{N}}$ such that
\begin{equation}
  \Gamma_i V_i \colon \mathbb{V}_{r_i} \to \mathbb{V}_{ r_{i+1} }\,,\quad
  \| Q^i \|_{r_i} < (q_i)^{-1}\,,\quad 
  |{Q^i}_{0,0}| > q_i\,,\quad 
  \| V_i \|_{r_i} < \epsilon_{\mu_i}^2
\end{equation}
where $r_i \equiv r-\sum_{j=1}^{i-1} \mu_j$. Moreover they must satisfy:
\begin{enumerate}[label=(\alph*)]
\item \label{convCond_i} $\epsilon_i\, =\,  q^3_i {\mu_i}^{2\tau+3}/(2\,C)$
                 as we computed in section \ref{sec:norms};
\item \label{c1} $ q_i\,\geq \tilde{C}/(C r_i^2) $, 
                 coherently with equation~\eqref{eq:convergenceCondition}; 
\item \label{c2} $0 < \mu_i < \tfrac{r_i}{3} $,
                 as required by Proposition~\ref{prop:kam};
\item \label{c3} $\sum_{i=1}^\infty \mu_i < r$,
                 to ensure that
                 \[
                   r_\infty=r-\sum_{j=1}^{\infty} \mu_j >0
                 \]
                 and so that the operator $\mathcal{H}_\infty$ 
		 is well defined on $\mathbb{V}_{r_\infty}$; 
\item \label{c4} $\lim_{i\to\infty} \epsilon_i = 0$, so 
                 we can conclude that $V_\infty \dfnd \lim_{i\to\infty}V_i = 0$;
\item \label{c5} $ 0 < q_\infty <q_i < 1,\quad \forall i \in \mathbb{N}$
                 as required by Proposition \ref{prop:bounds_rb};
\end{enumerate}
We choose 
\begin{equation}
  \epsilon_i\, =\, \frac{\epsilon_0}{(i+1)^{2(2\tau + 3) } },\qquad 
  \mu_i\, =\,\frac{ 1 }{ (1+i)^2 }\bigg( \frac{ 2\,C\, \epsilon_0 }{q_i^3} \bigg)^{ \tfrac{1}{2\tau+3} }
\end{equation}
so that condition \ref{convCond_i} is satisfied. Also condition \ref{c4} 
is evidently satisfied. Now we compute
\begin{equation}
  \sum_{i=0}^\infty \mu_i\,= 
    \,\sum_{i=0}^\infty \frac{1}{(1+i)^2}\bigg( \frac{2\,C\,\epsilon_0}{  q_i^3 } \bigg)^{ \tfrac{1}{2\tau+3} }\,\leq
      \, \bigg( \frac{ 2\,C\, \epsilon_0 }{ q_\infty^3 } \bigg)^{ \tfrac{1}{2\tau+3} }\,\sum_{i=1}^\infty \frac{1}{ i^2}
        \,\leq\, \frac{\pi^2}{6}\bigg( \frac{ 2\,C\, \epsilon_0 }{ q_\infty^3} \bigg)^{ \tfrac{1}{2\tau+3} }
\end{equation}
Both conditions \ref{c2} and \ref{c3} are satisfied by imposing 
$\sum_{i=0}^\infty \mu_i<r/3$ and, by the result above, we get
\begin{equation}
  \label{eq:eps0cond1}
    \bigg( \frac{ 2\,C\, \epsilon_0 }{ q_\infty^3} \bigg)^{ \frac{1}{2\tau+3} } \leq\, \frac{ 2\,r }{\pi^2}
\end{equation}
Then for the sequence $\{q_i\}_{i\in\mathbb{N}}$ we make the ansatz
\begin{equation}
  q_i = q_{i-1}\Big( 1-\frac{1}{(i+1)^2} \Big) = q_0 \prod_{j=2}^{i+1}\Big( 1-\frac{1}{j^2} \Big)
\end{equation}
By taking the logarithm of both sides, and using that \mbox{$\log(1-x)\geq -\log(4)x$} for 
\mbox{$x\in[\tfrac{1}{2},1]$}, we get
\begin{multline}
  \nonumber
  \log( q_i ) = \log(q_0) + \sum_{j=2}^{i+1} \log\Big( 1 - \frac{1}{j^2} \Big) \geq \\
    \geq \log(q_0) - \sum_{j=2}^{i+1} \frac{\log(4)}{j^2} 
      \geq \log(q_0) -\sum_{j=1}^\infty \frac{\log(4)}{j^2} = \log( q_0 \, 2^{-\pi^2/3} )
\end{multline}
We set $q_\infty = q_0 \, 2^{-\pi^2/3}$ and $q_0 < 1$ so condition \ref{c5} is satisfied.
If we plug this value for $q_\infty$ into equation~\eqref{eq:eps0cond1} we get
\begin{equation}
  \label{eq:eps0cond2}
    \epsilon_0\,\leq\,\frac{q_0^3}{C}\,\bigg(\frac{r}{\pi^2}\bigg)^{2\tau+3} 2^{ (4+2\tau-\pi^2)(2\tau+3)/(1-\pi^2) }
\end{equation}

Finally, we rewrite condition \ref{c1} as $r^2 \geq \tilde{C}/(qC)$ and, being 
$1/q \geq 1 $, we get a lower bound on $r$,
\begin{equation}
  r \geq \sqrt{ \tilde{C}/C}
\end{equation}

\end{proof}

\section{Conclusions}
\label{sec:conclusions}

So, in this paper we have provided an algorithm to perform 
perturbation theory for a Hamiltonian system on a Lie algebra. 
We assume to have a flow (unperturbed) that preserves some Lie 
subalgebra of the 
Lie algebra. When a perturbation is added to the 
given flow, the subalgebra is not preserved anymore. However, it
is possible to conjugate the perturbed flow to a new one, that 
preserves the same subalgebra of the unperturbed system,
up to terms quadratic in the perturbation.

while extending its scope beyond the original one.
Variants of the original theorem of Kolmogorov
\cite{kolmogorov_preservation_1954} have already been
proposed: for classical systems without action-angle 
coordinates \cite{de_la_llave_kam_2005}, for classical system
with degeneracy in the Hamiltonian\footnote{
  This means that the hessian of the Hamiltonian with
  respect to the action variables is not of maximal
  rank
} \cite{arnold_small_1963}, \cite{pinzari_properly-degenerate_2010},
or in theorin the volume preserving maps
and flows \cite{li_persistence_2002};
all of these cases are encompassed in our formula.
Nevertheless we can do more, and apply our method to 
non-canonical Poisson systems, which are gaining increasing 
importance in physics, since some pioneering works in the 1980ies
\cite{morrison_poisson_1982}, \cite{littlejohn_hamiltonian_1982}.

We have applied our theorem to the simple example of a 
non-autonomous symmetric Top. The Top is a non-canonical
Hamiltonian system with a degenerate bracket \eqref{eq:bracket},
which is not written in canonical coordinates. However, 
with a change of variables we can reduce it to a canonical 
form. When the moments of inertia have a prescribed 
time-dependence, the system becomes non-autonomous and 
is described by another angle variable (if it is 
periodically time-dependent).
One novelty with respect to classical mechanics is that the 
phase space has not the structure of a cotangent bundle.
We have shown that our formula can be iteratively
applied, to prove a KAM theorem for this dynamical
system. 

While on the one side we have shown that our method
fits in a typical KAM scheme, even if the system under
consideration fails the hypothesis of non-degeneracy,
we have not used many potentialities of our method. For
instance, we have introduced a set of canonical coordinates:
it would be 
interesting to reconsider the problem in the coordinates 
$(\coVec{M},t)$: this can be done with
our method, after a proper choice of the subalgebra $\mathbb{B}$ 
and of the operators $\mathcal{R}$ and $\Gamma$.
However, we think that the most interesting development 
would be to
write an iteration mechanism that works on any Lie algebra
to provide an algebraic KAM theorem.

\bibliography{bib_rigid_body}
\bibliographystyle{plain}

\appendix

\section{Proof of Proposition \ref{prop:bounds_rb}}
%\begin{proof}[Proof of Proposition \ref{prop:bounds_rb}]
\label{app:proof2}

By definition,
\begin{equation}
  (\Gamma W) Z\, =\, \lie{\mathcal{G}_s W} Z-(\mathcal{A}W)\,\partial_x Z-\lie{x\mathcal{G}_s Q \mathcal{A}W} Z - \lie{x\mathcal{G}_s Q \partial_\theta \mathcal{G}_s \mathcal{P}_0 W} Z
  \label{eq:theBeast}
\end{equation}
We will study each of the four terms on the r.h.s. separately. 
About the first one,
\begin{equation}
\footnotesize{
  \begin{aligned}
  \nonumber
  \|\lie{\mathcal{G}_s W} Z\|_{r -\delta -d}\,=\\
    =\,\sum_{L,M \in\mathbb{Z} } &\frac{e^{(r-\delta-d)(|L|+|M|)}}{\rho} \bigg| \sum_{l,m\in\mathbb{Z}_0} 
       \frac{ W_{l,m,1}\,(M-m)\,Z_{L-l,M-m}\,-\, m \,(\mathcal{P}_{\leq 1} W_{l,m})\,\partial_x Z_{L-l,M-m} }{ x_0 \Delta\rho m + l } \bigg|_{r-d} \leq\\[1ex]
  \nonumber
  \leq &\, \sum_{L,M \in\mathbb{Z}; l,m\in\mathbb{Z}_0}\, 
    e^{(r-\delta-d)(|L-l|+|M-m|)} \,e^{(r-\delta-d)(|l|+|m|)} 
     \frac{(|m|+|l|)^{\tau}}{\gamma\rho }\,\times \\[1ex]
  &\,\times \bigg( |W_{l,m,1}|_{r-\delta-d}\, |M-m|\,|Z_{L-l,M-m}|_{r-\delta-d}%\,+\\
       \,-\,|m|\,\big|\mathcal{P}_{\leq 1} W_{l,m})\big|_{r-\delta-d}\,\big|\partial_x Z_{L-l,M-m}\big| \bigg) \bigg|_{r-\delta-d}
       \leq \\[1ex]
  \nonumber
  \leq &\, \frac{1}{\rho\gamma}\sum_{L,M,l,m\in\mathbb{Z}} e^{r (|l|+|m|)}\,%\times \\[1ex]\,\times 
  \bigg( \frac{1}{ed}\, \Big(\frac{\tau}{e(d+\delta)} \Big)^{\tau}
    e^{(r-\delta)(|L-l|+|M-m|)}\, \frac{|W_{l,m}|_{r}}{d+\delta}\,|Z_{L-l,M-m}|_{r-\delta-d}\,+\\[1.5ex] 
  \nonumber
      &+\ \Big(\frac{\tau+1}{e(d+\delta)} \Big)^{\tau+1}\, e^{(r-\delta-d)(|L-l|+|M-m|)} 
        \,|W_{l,m}|_{r-\delta-d}\, \frac{|Z_{L-l,M-m}|_{r-\delta}}{d} \bigg) \leq
  \\[1.5ex]
  \end{aligned}
}
\end{equation}
In going from the 4th to the 5th line we used the condition \ref{1}. Now
\begin{equation}
  \label{eq:noUseLabel}
  \|\lie{\mathcal{G}_s W} Z\|_{r -\delta -d}\,
  \leq\, \frac{\tau^\tau + (\tau+1)^{(\tau+1)}}{\gamma\rho e^{\tau+1} (d+\delta)^{\tau+1} d}\,\|W\|_{r}\,\|Z\|_{r-\delta} 
  \equiv \frac{C_1}{(d+\delta)^{\tau+1} d}\,\|W\|_{r}\,\|Z\|_{r-\delta} 
\end{equation}
and in the last passage we introduced a constant $C_1$ for conciseness.
Next we consider
\begin{equation}
\footnotesize{
  \begin{aligned}
  \mathcal{A}W\,&=\,\frac{1}{Q_{0,0}}\bigg(\,\mathcal{P}_0\,\sum_{n\geq 0}\,(n+1)\, W_{0,0,n+1}\, x^n%\\[1ex]
     -\,\oint\sum_{ L,M\in\mathbb{Z} ; l,m\in\mathbb{Z}_0 }\,
    \frac{me^{iLt+iM\theta}}{x_0\rho\Delta m + l}\,Q_{L-l,M-m}\,W_{l,m,0}\bigg) \,=\\[1ex]
    &=\,\frac{W_{0,0,1}}{Q_{0,0}}\,-\,\frac{1}{Q_{0,0}}\sum_{l,m\in\mathbb{Z}_0}
      \frac{m}{x_0\rho\Delta m + l}\,Q_{-l,-m}\,W_{l,m,0}
  \end{aligned}
}
  \nonumber
\end{equation}
So that
\begin{equation}
  \notag
  \begin{aligned}
  |\mathcal{A}W|\, &\leq\ \Big|\frac{W_{0,0,1}}{Q_{0,0}}\Big|+\sum_{l,m\in\mathbb{Z}_0}
      \Big|\frac{m Q_{-l,-m} W_{l,m,0} }{Q_{0,0}( x_0\rho\Delta m + l) }\Big|\leq\\[1ex]
  &\leq\ \frac{1}{q}\, \|\partial_x W_{0,0}\|_{r-\delta-d}\, +\, 
    \frac{1}{q\gamma}\,\sum_{l,m\in\mathbb{Z}_0}\,
      (|m| +| l|)^{(\tau+1)}\,\|\mathcal{P}_0 W\|_{r}\, e^{-r(|m|+|l|)}\,|Q_{-l,-m}| 
  \end{aligned}
\end{equation}
where we used $ |\mathcal{P}_0 W_{l,m} |_{r} \leq \|\mathcal{P}_0 W\|_{r} e^{-r(|l|+|m|)} $;
also, in passing from the first to the second line we employed hypothesis \ref{2}.
Continuing:
\begin{equation}
  \nonumber
  \begin{aligned}
  |\mathcal{A}W|\ &
    \leq\ \frac{\|W\|_{r}}{q(d+\delta)} 
    +\,\frac{(\tau+1)^{\tau+1}}{q \gamma (e r)^{\tau+1}}
     \|W\|_{r}\sum_{l,m\in\mathbb{Z}_0}|Q_{l,m}|
       \leq \\&
  \leq\ \frac{\|W\|_{r}}{q (d+\delta)}\, 
    +\,\frac{(\tau+1)^{\tau+1}}{q^2 \gamma (e r)^{\tau+1}}
     \|W\|_r 
  \leq\ \frac{ C_2 \|W\|_{r} }{ q^2 (d+\delta)^{\tau+1} }
  \end{aligned}
\end{equation}
where in the last passage we used 
\mbox{$\delta+d\leq r$} so that \mbox{$r^{-1}\leq(d+\delta)^{-1}$},
and $C_2$ is a constant. 
Then for the second term of equation~\eqref{eq:theBeast} we have
\begin{equation}
  \label{eq:a_bound}
  \|(\mathcal{A}W)\partial_x Z\|_{r -\delta - d}\,\leq\, \sum_{l,m\in\mathbb{Z}} e^{r(|l|+|m|)}\, |aW|\, \frac{1}{d}\, |Z_{l,m}|_{r-\delta}\,%\leq \\[1ex]
      \leq\, \frac{C_2 \|W\|_{r}\|Z\|_{r-\delta}}{q^2\,d\,(d+\delta)^{\tau+1}}
\end{equation}

The third term of equation~\eqref{eq:theBeast} reads 
\begin{equation}
%  \nonumber
  \lie{x\mathcal{G}_sQ\mathcal{A}W}Z\, =\, \frac{1}{\rho}\sum_{L,M\in\mathbb{Z}}\,e^{iLt+iM\theta}(\mathcal{A}W)\,\sum_{l,m\in\mathbb{Z}_0}
  \Big(\footnotesize{
      \frac{Q_{l,m}\,(M-m)\,Z_{L-l,M-m}\,-\,x\,m\,Q_{l,m}\,\partial_x Z_{L-l,M-m}}{x_0\rho\Delta m + l}
  }\Big) 
\nonumber
\end{equation}
so that 
\begin{equation}
\footnotesize{
  \notag \nonumber
  \begin{aligned}
  \big\| \lie{x\mathcal{G}_sQ & \mathcal{A}W}Z \big\|_{r-\delta -d}\, 
    \leq \\[1ex]
  \leq& \,\frac{|\mathcal{A}W|}{\rho}\,\sum_{L,M\in\mathbb{Z}}\,e^{(r-\delta-d)(|L|+|M|)}\, 
      \bigg|\sum_{l,m\in\mathbb{Z}_0} \frac{ Q_{l,m}\,(M-m)\,Z_{L-l,M-m}\,
        -\,x\,m\,Q_{l,m}\,\partial_x Z_{L-l,M-m}}{x_0\rho\Delta m + l} \bigg|_{r-\delta-d}\,\leq \\[1ex]
  \leq& \,\frac{|\mathcal{A}W|}{\rho\gamma}\,\sum_{L,M,l,m\in\mathbb{Z}}\, e^{(r-\delta-d)(|L-l|+|M-m|)}\,
    \,e^{(r-\delta-d)(|l|+|m|)}|Q_{l,m}| \times\\[1ex]
  &\,\times\bigg(|M-m|\, \frac{(|m|+|l|)^\tau}{\gamma} \,|Z_{L-l,M-m}|_{r-\delta-d}\,
  +\,|x|_{r-\delta-d}\,\,\big|\partial_x Z_{L-l,M-m}\big|_{r-\delta-d}\, \frac{(|m|+|l|)^{\tau+1}}{\gamma}\,\bigg)\leq\\[1ex] 
  \leq & \frac{|\mathcal{A}W|}{\rho\gamma}\,\bigg(\frac{1}{ed}\,\|Z\|_{r-\delta}\,\Big(\frac{\tau}{e(d+\delta)}\Big)^\tau
    \|Q\|_r\,
  +\,\frac{|r|}{d}\,\Big(\frac{\tau+1}{e(d+\delta)}\Big)^{\tau+1}\|Z\|_{r-\delta}\|Q\|_r\bigg)\,\leq\, 
  \,\frac{C_3\|W\|_{r}\|Z\|_{r-\delta} }{q^3\,d\,(d+\delta)^{2\tau+2}}
  \end{aligned}
}
\end{equation}
where $C_3$ is another constant.

Finally, the fourth term of equation~\eqref{eq:theBeast} is
\begin{equation}
\footnotesize{
  \notag
  \begin{aligned}
  \lie{x & \mathcal{G}_sQ\partial_\theta\mathcal{G}_s\mathcal{P}_0W}Z\,
    = \, \bigg\{ x \mathcal{G}_s \sum_{L,M,l,m \in \mathbb{Z}_0 } \, e^{iLt+iM\theta} 
    \frac{ m W_{l,m,0} Q_{L-l,M-m} }{ x_0\rho\Delta m + l }
    \,\sum_{l_1,m_1\in\mathbb{Z}} Z_{l_1,m_1} e^{il_1t+im_1\theta} \bigg\}\, =\\[1ex]
    =&\,\sum_{L_1,M_1\in\mathbb{Z};L,M,l,m\in\mathbb{Z}_0}\,e^{iL_1 t+iM_1\theta}
    \frac{ 
      m\,W_{l,m}\, Q_{L-l,M-m}\,(M_1-M)\,Z_{L_1-L,M_1-M}\,
      -\,x\,M\,m\,W_{l,m}\,Q_{L-l,M-m}\,\partial_xZ_{L_1-L,M_1-M}
    }{ 
      \rho (x_0\rho\Delta M + L) (x_0\rho\Delta m + l) 
    }
  \end{aligned}
}
\end{equation}
and so
\begin{multline}
  \|\lie{x\mathcal{G}_s  Q \partial_\theta \mathcal{G}_s \mathcal{P}_0 W} Z \|_{r-\delta-d} \leq \\
  \leq 
  \sum_{L_1,M_1,L,M,l,m\in\mathbb{Z}}\bigg( e^{(r-\delta-d)(|L_1 - L|+|M_1-M|)}\, e^{(r-\delta-d)(|L|+|M|)}\,|W_{l,m,0}|\,|Q_{L-l,M-m}|\,
    \frac{(|m|+|l|)^{\tau+1}}{\rho\gamma^2}  \times \\[1ex] \times \Big( |M_1-M|\, 
      |Z_{L_1-L,M_1-M}| \,(|M|+|L|)^\tau\, +\, |x|(|M|+|L|)^{\tau+1} |\partial_x Z_{L_1-L,M_1-M}|_{r-\delta-d} \bigg) \,\leq\\[1ex]
  \leq \frac{\|Z\|_{r-\delta}}{q\rho\gamma^2 d e^{\tau+1} }\,  \|W\|_{r}
    \bigg( \bigg( \frac{2\tau}{d+\delta}\bigg)^\tau \bigg(\frac{2(\tau+1)}{e(d+\delta)}\bigg)^{\tau+1}
      + |r| \bigg(\frac{2(\tau+1)}{d+\delta}\bigg)^{\tau+1}\bigg(\frac{2(\tau+1)}{e(d+\delta)}\bigg)^{\tau+1} \bigg) \equiv \\[1ex]
  \,\equiv\, \frac{ C_4 \|W\|_{r} \|Z\|_{r-\delta} }{q\,d\,(d+\delta)^{2\tau +1 } }
  \nonumber
\end{multline}
with a fourth constant $C_4$. By defining
\begin{equation}
  C\, = \,(C_1\,q^3 + C_2\,q)(d+\delta)^{\tau+1} + C_3 + (d+\delta) q^2 C_4
\end{equation}
we end up with the thesis.

To prove the second and third inequalities,
we start by observing that
\begin{eqnarray} 
%  \nonumber
  &\|\mathcal{N}_s V \|_{r-\mu} = \|\chi\mathcal{P}_0 V + \mathcal{P}_1 V\|_{r-\mu}\leq \|V \|_{r-\mu}\\
%  \nonumber
  &\|\mathcal{R}_s V \|_{r-\mu} = \|\oint\mathcal{P}_0 V + \mathcal{P}_{2} V\|_{r-\mu}\leq \|V \|_{r-\mu} 
\end{eqnarray}

Then we consider
\begin{equation}
  \| \mathcal{K}V \|_{r-\mu}\, \leq\, | x_0 \Delta \rho \mathcal{A} V |\, +\,
    \big\|\lie{\tfrac{1}{2}Qx^2} x\mathcal{G}_s( \mathcal{P}_0 \partial_x - Q \mathcal{A} - Q \partial_\theta \mathcal{G}_s \mathcal{P}_0 ) V \big\|_{r -\mu} 
\end{equation}
By~\eqref{eq:a_bound},
\begin{equation}
  | x_0 \Delta \rho \mathcal{A}V | \,\leq\, |x_0|\, \rho \Delta\, \big| \mathcal{A}V \big|\, \leq\, \rho\Delta\, \frac{C_2}{q^2\,\mu^{\tau +1 }}\, \|V\|_{r} 
\end{equation}

To the next term we apply equation~\eqref{eq:formula3} of 
Proposition~\ref{prop:normProp} with $\delta=d=\mu/2$,
\begin{multline}
  \nonumber
  \big\|\lie{\tfrac{1}{2}Qx^2} x\mathcal{G}_s( \mathcal{P}_0 \partial_x - Q \mathcal{A} - Q \partial_\theta \mathcal{G}_s \mathcal{P}_0 ) V \big\|_{r -\mu} \leq\\ 
    \leq \frac{4}{\rho e \mu^2} \big\|\tfrac{1}{2}Qx^2 \big\|_{r}\, 
      \big\|x\mathcal{G}_s( \mathcal{P}_0 \partial_x - Q \mathcal{A} - Q \partial_\theta \mathcal{G}_s \mathcal{P}_0 ) V \big\|_{r -\mu/2}
\end{multline}
We have
\begin{equation}
  \big\| \tfrac{1}{2}Qx^2 \big\|_r\,\leq\,\tfrac{1}{2}\|Q\|_r |x^2|_r \leq (2q)^{-1}
\end{equation}
and
\begin{equation}
  \nonumber
  \begin{aligned}
    \big\| x \mathcal{G}_s & \mathcal{P}_0 \partial_x V \big\|_{r - \mu/2 }\,=\sum_{l,m\in\mathbb{Z}_0}\Big| 
      \frac{ V_{lm1}x }{ x_0\rho\Delta m + l } \Big|_{ r-\tfrac{\mu}{2} }\, e^{ ( r-\mu/2 ) (|l|+|m|) }\leq \\[1ex]
        &\leq\,\sum_{ l,m\in\mathbb{Z}_0 } | V_{lm1}x |_{ r-\mu/2 }\,\big(|m| + |l|\big)^{\tau} \gamma^{-1}\,e^{ (r-\mu/2) (|l|+|m|) }
          \,\leq\,\frac{1}{\gamma}\Big(\frac{2\tau}{e\mu}\Big)^{\tau+1} \|V\|_{r}
  \end{aligned}
\end{equation}
and also
\begin{equation}
  \nonumber
  \begin{aligned}
    \big\| x\mathcal{G}_s Q \mathcal{A} V \big\|_{r-\mu/2}\, &=\, 
      \big|\mathcal{A}V\big|\, |x|_{r-\mu/2}\, \big\| \mathcal{G}_s Q \big\|_{r-\mu/2}\,\leq \\[1ex] 
      &\leq \big| \mathcal{A}V \big|\sum_{l,m\in\mathbb{Z}_0} \gamma^{-1}(|l|+|m|)^\tau e^{(r-\mu/2)(|l|+|m|)}|Q_{l,m}|\,\leq \\[1ex]
      &\leq\,\frac{ C_2 }{ \gamma q^2 \mu^{\tau+1} }\,\|V\|_r\, \frac{1}{}\,\bigg(\frac{2\tau}{e\mu}\bigg)^\tau\|Q\|_r
      \leq\,\frac{ C_2 }{ \gamma q^3\mu^{\tau+1} }\,\bigg( \frac{ 2\tau }{ e\mu }\bigg)^\tau\|V\|_{r}
  \end{aligned}
\end{equation}
And finally
\begin{equation}
\footnotesize{
  \notag
  \begin{aligned}
%  \nonumber
  \big\| x\mathcal{G}_s & Q\partial_\theta \mathcal{G}_s \mathcal{P}_0 V \big\|_{r-\tfrac{\mu}{2}}\,%=\\[1ex]
    =\sum_{L,M\in\mathbb{Z}_0} \frac{ e^{(r-\tfrac{\mu}{2})(|L|+|M|)} }{| x_0\Delta\rho M + L |}\bigg| \sum_{l,m\in\mathbb{Z}_0} 
      \frac{m Q_{L-l,M-m} v_{l,m,0} }{ x_0\Delta\rho m + l } \bigg|\leq\\[1ex]
    &\leq\sum_{L,M,l,m\in\mathbb{Z}} e^{(r-\tfrac{\mu}{4} )(|L|+|M|)} e^{-\tfrac{\mu}{4}(|L|+|M|)}
       \frac{ (|M| + |L|)^{\tau} }{\gamma}%\times\\[1ex] \times 
  \frac{(|m|+|l|)^\tau  }{ \gamma } |Q_{L-l,M-m}| |V_{l,m,0} |\,\leq\\[1ex]
  &\leq \, \Big( \frac{4\tau}{e\mu}\Big)^\tau \sum_{L,M,l,m\in\mathbb{Z}} 
    e^{ r (|L|+|M|)-\tfrac{\mu}{4}( |L-l|+|M-m| )-\tfrac{\mu}{4}( |l|+|m| ) } \frac{(|m|+|l|)^\tau  }{\gamma^2} |Q_{L-l,M-m}| |V_{l,m,0} | \leq\\[1ex]
    &\leq \Big( \frac{4\tau}{e\mu}\Big)^\tau \Big( \frac{4(\tau+1) }{e\mu}\Big)^{\tau+1}\frac{ \|V\|_{r} }{\gamma^2 q} 
  \end{aligned}
}
\end{equation}
We can conclude that 
\begin{equation}
  \| \mathcal{K}V \|_{r-\mu} \leq \frac{\tilde{C} \|V\|_{r} }{q^3 \mu^{2\tau+3}} \\
\end{equation}
and so
\begin{equation}
  \| \mathcal{N}V \|_{r-\mu} \,\leq\, \| \mathcal{N}_s V\|_{r-\mu}\, +\, \| \mathcal{K}V \|_{r-\mu}\,
    \leq\, \|V\|_r \, +\, \frac{ \tilde{C}_1 \|V\|_{r} }{ q^3 \mu^{2\tau+3} } 
      \,\equiv\, \frac{ \tilde{C} \|V\|_r }{ q^3\, \mu^{2\tau+3} }
\end{equation}
and analogously
\begin{equation}
  %\nonumber
  \| \mathcal{R}V \|_{r-\mu} \,\leq\, \| \mathcal{R}_s V\|_{r-\mu}\, +\, \| \mathcal{K}V \|_{r-\mu}\,
    \leq\, \|V\|_{r} \, +\, \frac{ \tilde{C}_1 \|V\|_r }{ q^3 \mu^{2\tau+3} } 
      \,\equiv\, \frac{ \tilde{C} \|V\|_{r} }{q^3\, \mu^{2\tau+3} }
\end{equation}
where \mbox{$\tilde{C}_1,\tilde{C}$} are constants depending on $\mu,\tau,\gamma,q,\rho,\Delta,e$.
This concludes the proof.

\end{document}